\documentclass[mathpazo]{csam}

\usepackage{latexsym}
\usepackage{algpseudocode}

\usepackage{graphicx, color}
\usepackage{multirow}
\usepackage{subfigure}

\usepackage{amsfonts}
\usepackage{amsmath,amssymb}
\usepackage{amsthm}
\usepackage{mathtools}
\usepackage{enumitem}

\usepackage{algorithm}
\usepackage{algorithmicx}
\usepackage{algpseudocode}

\usepackage{float}
\usepackage{epsfig}
\usepackage{array,booktabs}
\usepackage{tabularx}
\usepackage[pagewise]{lineno}  % \linenumbers
\usepackage{MnSymbol}
\providecommand{\U}[1]{\protect\rule{.1in}{.1in}}

\usepackage{longtable}

\usepackage{amsopn}

\newcommand{\RomanNumeralCaps}[1]
    {\MakeUppercase{\romannumeral #1}}

\newtheorem{assumption}{Assumption}

\newdimen\slantmathcorr
\def\oversl#1{%assuming that mathslant=0.25
\setbox0=\hbox{$#1$}
\slantmathcorr=\wd0
\hskip 0.2\slantmathcorr \overline{\hbox to 0.8\wd0{%
\vphantom{\hbox{$#1$}}}}
\hskip-\wd0\hbox{$#1$}
}

\def\undersl#1{%assuming that mathslant=0.25
\setbox0=\hbox{$#1$}
\slantmathcorr=\wd0
\underline{\hbox to 0.8\wd0{%
\vphantom{\hbox{$#1$}}}}
\hskip-0.8\wd0\hbox{$#1$}
}

\begin{document}

%%%%% title : short title may not be used but TITLE is required.
% \title{TITLE}
% \title[short title]{TITLE}
\title{Asymptotic expansion regularization for inverse source problems in two-dimensional singularly perturbed nonlinear parabolic PDEs}

%%%%% author(s) :
% single author:
% \author[name in running head]{AUTHOR\corrauth}
% [name in running head] is NOT OPTIONAL, it is a MUST.
% Use \corrauth to indicate the corresponding author.
% Use \email to provide email address of author.
% \footnote and \thanks are not used in the heading section.
% Another acknowlegments/support of grants, state in Acknowledgments section
% \section*{Acknowledgments}
%\author[O.~Author]{Only Author\corrauth}
%\address{School of Mathematical Sciences, Beijing Normal University,
%Beijing 100875, P.R. China}
%\email{{\tt author@email} (O.~Author)}

% multiple authors:
% Note the use of \affil and \affilnum to link names and addresses.
% The author for correspondence is marked by \corrauth.
% use \emails to provide email addresses of authors
% e.g. below example has 3 authors, first author is also the corresponding
%      author, author 1 and 3 having the same address.
% \author[Z. Zhang et~al.]{Zhengru Zhang\affil{1}\comma\corrauth,
%       Author Chan\affil{2}~and Author Zhao\affil{1}}
% \address{\affilnum{1}\ School of Mathematical Sciences,
%          Beijing Normal University,
%          Beijing 100875, P.R. China. \\
%           \affilnum{2}\ Department of Mathematics,
%           Hong Kong Baptist University, Hong Kong SAR.}
% \emails{{\tt zhang@email} (Z.~Zhang), {\tt chan@email} (A.~Chan),
%          {\tt zhao@email} (A.~Zhao)}
% \footnote and \thanks are not used in the heading section.
% Another acknowlegments/support of grants, state in Acknowledgments section
% \section*{Acknowledgments}

\author[D. Chaikovskii et~al.]{Dmitrii Chaikovskii\affil{1},       Aleksei Liubavin\affil{2}~and Ye Zhang\affil{1}\comma\affil{3}\comma\corrauth}
 \address{\affilnum{1}\ Shenzhen MSU-BIT University,  Shenzhen 518172, China. \\
            \affilnum{2}\ School of Mathematical Sciences, East China Normal University, Shanghai, 200241, China. \\
           \affilnum{3}\ School of Mathematics and Statistics, Beijing Institute of Technology, Beijing, 100081, China.}
           
\emails{{\tt dmitriich@smbu.edu.cn} (D.~Chaikovskii), {\tt la1992@mail.ru} (A.~Liubavin),   {\tt ye.zhang@smbu.edu.cn} (Y.~Zhang)}

%%%%% Begin Abstract %%%%%%%%%%%
\begin{abstract}
In this paper, we develop an asymptotic expansion-regularization (AER) method for inverse source problems in two-dimensional nonlinear and nonstationary singularly perturbed partial differential equations (PDEs). The key idea of this approach is the use of the asymptotic-expansion theory, which allows us to determine the conditions for the existence and uniqueness of a solution to a given PDE with a sharp transition layer. As a by-product, we derive a simpler link equation between the source function and first-order asymptotic approximation of the measurable quantities, and based on that equation we propose an efficient inversion algorithm, AER, for inverse source problems. We prove that this simplification will not decrease the accuracy of the inversion result, especially for inverse problems with noisy data. Various numerical examples are provided to demonstrate the efficiency of our new approach.
\end{abstract}
%%%%% end %%%%%%%%%%%

%%%%% AMS/PACs/Keywords %%%%%%%%%%%
%\pac{}
\ams{65M32, 35C20, 35G31%The information of the AMS subject classification can be found in http://mathscinet.ams.org/msc/msc2010.html
}

%%%% maketitle %%%%%
\maketitle

%%%% Start %%%%%%
\section{Introduction}\label{sec:Introduction}
In this paper, we develop an asymptotic expansion-\hspace{0cm}regularization method for two-\hspace{0cm}dimensional inverse source problems which arise from time-dependent singularly perturbed partial differential equations (PDEs).  To illustrate our ideas, we take the following inverse problem as an example:\\
(\textbf{IP}): Given noisy data $\{u^\delta(x,y,t_0), u_x^\delta(x,y,t_0), u_y^\delta(x,y,t_0)\}$ of $\{u(x,y,t),$ $ u_x(x,y,t),$ $ u_y(x,y,t)\}$ at the $n\cdot m$ location points $\{x_i, y_j\}^{n,m}_{i,j=0}$ and at the time point $t_0$, find the source function $f(x,y)$ such that $(u,f)$ satisfies the dimensionless nonlinear autowave model
\begin{align} \label{mainproblem}
\begin{cases}
\displaystyle \mu  \Delta u-\frac{\partial u}{\partial t} =-u \left( k \frac{\partial u}{\partial x}+ \frac{\partial u}{\partial y} \right)+f(x, y), \quad  x \in \mathbb{R}, \quad  y \in  (-a, a) , \quad  t\in (0, T], \\
u(x, y, t, \mu)=u(x+L, y, t, \mu), \quad x\in \mathbb{R}, \quad y\in[-a, a],\quad t\in[0, T], \\
u(x, -a, t, \mu)=u^{-a}(x), \quad u(x, a, t, \mu)=u^{a}(x), \quad x\in \mathbb{R}, \quad t\in[0, T], \\
u(x, y, 0, \mu)=u_{init}(x, y, \mu), \quad  x\in \mathbb{R}, \quad y\in[-a, a],
\end{cases}
\end{align}
where $ u(x,y,t)$ represents the temperature or oil saturation, $\mu \ll 1$ denotes kinematic viscosity, the positive constant $k$ is the medium anisotropy coefficient, and $f(x,y)$ is the source function.
We assume that the function $f( x, y)$ is $L$-periodic in the variable $x$ and sufficiently smooth in the region $(x,y):\mathbb{R} \times \bar{\Omega}$ ($\Omega\equiv(-a, a)$), that  the functions $u^{-a}(x)$, $u^{a}(x)$ are $L$-periodic and sufficiently smooth in $x\in \mathbb{R}$, and that $u_{init}(x, y, \mu)$  is a sufficiently smooth function in $(x,y):\mathbb{R} \times \Omega$ and L-periodic in $x$, satisfying  $u_{init}(x, -a, \mu)=u^{-a}(x), u_{init}(x, a, \mu)=u^{a}(x)$. In this paper, we focus on the speed, location, and width of the border between two regions - the region with a small dimensionless value  $u$ and the region with its high value. The domain of the function describing the moving front contains a subdomain in which the function has a large gradient. Interest in front-type solutions is associated with combustion problems \cite{Liberman2003} or nonlinear acoustic waves \cite{Rudenko2017}.

Note that the inverse source problem (\textbf{IP}) is ill-posed (see \cite{Isakov1990}); we should therefore employ the regularization methods to obtain a meaningful approximate source function. Within the framework of Tikhonov regularization, the (\textbf{IP}) can be converted to the following PDE-constrained optimization problem:
\begin{align} \label{LSM0}
\min_{f} \sum^n_{i=0} \sum^m_{j=0} \left\{ \left[u(x_i,y_j,t_0) - u^\delta(x_i,y_j,t_0) \right]^2 + \left[\frac{\partial u}{\partial x}(x_i,y_j,t_0) - \frac{\partial u^\delta}{\partial x}(x_i,y_j,t_0) \right]^2  \right. \nonumber \\ \left. + \left[\frac{\partial u}{\partial y}(x_i,y_j,t_0) - \frac{\partial u^\delta}{\partial y}(x_i,y_j,t_0) \right]^2  \right\} + \varepsilon \mathcal{R}(f),
\end{align}
where $u$ solves the nonlinear PDE \eqref{mainproblem} with a given $f$, $\mathcal{R}(f)$ denotes the regularization term, and $\varepsilon>0$ is the regularization parameter.

Although the conventional formulation \eqref{LSM0} for (\textbf{IP}) is straightforward, the numerical realization is difficult in many applications for the following three reasons: (1) The regularization term $\mathcal{R}(f)$ reflects the a priori information about the source function, which is hard to obtain in practice. (2) The regularization parameter $\varepsilon$ also needs to be carefully chosen, which is not an easy task. For example, the most popular a posteriori regularization-selection methods, e.g. the Morozov discrepancy principle, require repeated solving of the nonlinear high-order PDE \eqref{mainproblem}, which is time-consuming. (3) The optimization problem is highly non-convex, and hence, in general, there is no efficient algorithm for solving \eqref{LSM0} due to the existence of many local minima. To overcome these difficulties, in this paper, by using the asymptotic expansion method, we develop a link equation between the unknown source function and observation quantities, which avoids the resolution of high-order PDE \eqref{LSM0}.

In the literature, there are many examples of applications of the reaction–\hspace{0cm}diffusion–\hspace{0cm}advection equations with small parameters in problems of biology  \cite{PattersonWagner2012,DoOwida2011,BodnarSequeira2008,HidalgoTello2014}, physics \cite{BerrymanHolland1978}, and chromatography \cite{zhang2016regularization,zhang2016,LinZhang2018,ChengLin2018}, as well as in industrial problems \cite{PhysRevE.70.026307,2004evapormetals,2014evolutiondisp}. In addition, there is extensive literature devoted to the formation and dynamics of the combustion front using the reaction–diffusion–advection equations \cite{Vladimirova2006FlameCW,Calder2007CapturingTF,Libermanarticle}. Due to the depletion of fields with readily available oil reserves, the world's leading oil and gas companies have been forced to move to the development of fields with hard-to-recover reserves. It has thus become necessary to extract hydrocarbons from complex reservoirs with low porosity and permeability, and high-viscosity oils; in turn, the development and improvement of enhanced oil-recovery methods have become matters of some urgency. One such class of methods is that of thermochemical methods, which are based on the ability of reservoir oil to react with oxygen injected into the reservoir, accompanied by the release of a large amount of heat, i.e. in-situ combustion. According to the experimental data \cite{Burger1972,Burger2013}, under certain conditions, sharp transition layers appear in the reservoir in terms of temperature, oxidant concentration, and oil saturation. The formation and dynamics of such fronts can be modeled using reaction–diffusion–advection problems. Singularly perturbed problems are very well suited to the modeling of autowave problems, since the sharp transient layer formed by solving such problems creates a boundary between two regions in different states. Thus, in the article \cite{b04}, the authors describe the recent construction of an autowave model predicting the growth of the city of Shanghai in the near future. Problems with a solution in the form of a front on a segment are considered in \cite{NEFEDOV201390,b6}, and the motion of a two-dimensional front is studied in \cite{LevashovaTwoDim2021,AntLevNef18,AntVolLev17,Nefedov2020Cubic}.

The numerical methods for singularly perturbed problems can involve difficulties in making certain choices, such as the choice of the initial conditions lying in the area of influence of the solution with an inner transition layer, and the selection of adequate grids for the implementation of difference schemes. An analytical study of the solution is an effective means of overcoming these difficulties. The asymptotic methods used in this work, in particular Vasil'eva's algorithm \cite{ButuzovVasileva1970} and the asymptotic method of differential inequalities \cite{Nefedovarticle2000,Nefedov2008}, allow us to determine, up to a small parameter, the position of the transition layer and the equation of its motion \cite{Lukyanenkoarticle2016,Lukyanenko2017,Lukyanenkoarticle2017}, and to substantiate the existence of a solution of the considered type and thereby confirm the reliability of numerical calculations.

Of course, there are various other numerical methods for solving the reaction–\hspace{0cm}diffusion–advection equations, such as finite-difference methods, finite-element methods, and finite-volume methods. However, the asymptotic methods have their own merits. The asymptotic representation of the solution to the reaction–diffusion–advection equations makes it possible to describe the dynamics of a sharp transition layer. Therefore, we can estimate its width and also determine the shape of the front at each moment in time. It should be noted that the indicated asymptotic representation is rather simple, which is extremely important for obtaining an efficient algorithm for solving the inverse problem of finding the source function. Asymptotic expansion methods significantly accelerate the obtaining of approximate solutions with a required calculation accuracy, which leads to an increase in the efficiency of numerical calculations.

The main contribution of this paper is the development of a new methodology for efficiently solving inverse problems in singularly perturbed PDEs. To achieve this, we use the asymptotic method to reduce the original high-derivative PDE with a small parameter to a simpler equation of lower degree while obtaining sufficiently accurate results. It should be noted that a similar idea was used in \cite{LukyanenkoTwoDimInv2018} for a parameter-identification problem in a singularly perturbed reaction–diffusion equation, and in \cite{VolkovNevedofCoefInvProbl2020} for a one-dimensional coefficient inverse problem in a Burger-type equation.

The remainder of this paper is structured as follows. In Section \ref{constructionOfAsymptotic}, we perform asymptotic analysis for PDE \eqref{mainproblem}. The main results are presented in Section \ref{statementresults}, while Section  \ref{derivationAndProofs} provides technical proofs of the main theoretical results. In Section \ref{simulation}, some experiments for the forward and inverse problems are described. Finally, concluding remarks are made in Section \ref{Conclusion}.

\section{Asymptotic analysis for the forward problem}
\label{constructionOfAsymptotic}

In this section, we construct an asymptotic solution of PDE \eqref{mainproblem}, whose properties will be used in the subsequent sections. We focus on the solution to problem \eqref{mainproblem}, having the form of a moving front: at each moment in time at $ -a\leq y < h(x, t)$ the solution is close to the surface $\varphi^{(-)} (x, y)$, and at $h(x, t)< y\leq a$ it is close to the surface $\varphi^{(+)}(x, y)$, and sharply changes from the values of the surface $\varphi^{(-)} (x, y )$ to values of the surface $\varphi^{(+)}(x, y)$ in the neighborhood of the curve $y=h(x, t)$. In this case, the solution to problem \eqref{mainproblem} has an inner transition layer in the vicinity of this curve.

Above all, we list the assumptions under which the asymptotic solution exists:

\begin{assumption}
\label{A1}
$u^{-a}(x)<0$, $u^{a}(x)>0$ and $u^{a}(x)-u^{-a}(x)> 2\mu^2 $ for all $x \in \mathbb{R}$.
\end{assumption}

\begin{assumption}
\label{A2}
For any $x\in \mathbb{R}, y \in \bar{\Omega}$, the following holds:
\begin{align*}
\left(u^{-a}(x- k (a+y))\right)^2 > \frac{2}{k} \int_{x}^{x-k(a+y)} f\left(s,\frac{-x+ky+s}{k} \right) ds,   \\
\left(u^{a}(x + k(a -y))\right)^2 > \frac{2}{k} \int_x^{x+k(a -y )} f\left( s,\frac{-x+k y+s}{k} \right) ds.
\end{align*}

For convenience, for $(x,y) \in \mathbb{R}\times \bar{\Omega} \equiv S  $, we can rewrite Assumption $\ref{A2}$ as

\begin{align*}
\min_{x\in \mathbb{R}} \left( \left( u^{-a}(x)\right)^2\right) > - \frac{2}{k} \iint\limits_S \min \left(0,  f\left(x,y\right)\right) dx \, dy ,   \\
\min_{x\in \mathbb{R}} \left( \left( u^{a}(x)\right)^2\right) > \frac{2}{k} \iint\limits_S  \max \left(0,  f\left(x,y\right)\right)dx \, dy .
\end{align*}

\end{assumption}
\begin{assumption}\label{A3}
  For all $(x,t) \in \mathbb{R}\times \mathcal{T}$, $-a < {h_0}(x, t)<a$, where $h_0$ is the zero approximation of $h$ (see \eqref{curveexpansion}), and $ \displaystyle \max_{x\in \mathbb{R},\ t\in \bar{\mathcal{T}} } \left(  \frac{\partial {h_0} (x,t)}{\partial x}  \right) <1/k $.
\end{assumption}
\begin{assumption}
\label{A4}
$u_{init}(x,y)=U_{n-1}(x,y,0)+\mathcal{O}(\mu^n)$, where the asymptotic solution $U_{n-1}$ will be constructed later; see, e.g., Theorem \ref{MainThm}.
\end{assumption}

Under Assumptions \ref{A1} and \ref{A2}, the zero-order outer functions $\varphi^{(-)} (x,y)$ and $\varphi^{(+)}(x,y)$, which will be used for asymptotic construction (cf. \eqref{u0regu}), can be expressed explicitly from equations \eqref{zeroorderregularequation1} and \eqref{zeroorderregularequation2} in the following form:
\begin{align}\label{eq7}
\begin{split}
\varphi^{(-)} (x,y) = -\sqrt{\left(u^{-a}(x- k (a+y))\right)^2-\frac{2}{k} \int _x^{x- k (a+y)} f\left(s,\frac{-x+k y+s}{k}\right)ds},\\
\varphi^{(+)} (x,y) = \sqrt{\left(u^{a}(x + k(a -y))\right)^2-\frac{2}{k} \int _x^{x + k(a -y)} f\left(s,\frac{-x+k y+s}{k}\right)ds}.
\end{split}
\end{align}

Assumption \ref{A3} determines  the location of the inner transition layer within the specified region $\Omega$ for the variable $y$ and ensures that the transition-layer functions $Q_0^{(\mp)}(\xi,l,h_{0},t)$ are increasing (see \eqref{Q0equation}), while Assumption \ref{A4} means that at $t=0$ the transition layer has already been formed, and the initial function $u_{init}$ already has the transition layer in the vicinity of the curve $h_{0}^{*}:= h_0(x,0)$.

The curve $y=h(x, t)$ at each moment in time divides the region $\bar{\Omega}$ into two parts: $\bar{\Omega}^{(-)} =\{y:y\in[-a,h(x, t) ] \}$ and $\bar{\Omega}^{(+)}=\{y:y\in[h(x, t),a]\}.$

For a detailed description of the transition layer, in the vicinity of curve $ h (x, t)$, we proceed to the extended variable
\begin{equation} \label{xidefin}
\displaystyle \xi=\frac{r}{\mu},
\end{equation}
and to the local coordinates $ (l, r) $ using the relations
\begin{equation} \label{localcoord}
 x=l-r\sin\alpha, \quad y=h(l, t)+r\cos\alpha,
\end{equation}
where
\begin{equation} \label{sincosalpha}
\displaystyle \sin\alpha=\frac{h_{l}}{\sqrt{1+h_{l}^{2}}}, \quad \cos\alpha=\frac{1}{\sqrt{1+h_{l}^{2}}},
\end{equation}
$\alpha$ is the angle between axis $y$ and the normal to the curve $ y = h(x, t) $ drawn to the region $ y> h(x, t) $ for each $ t $, plotted counterclockwise, $ l $ is the $x$ coordinate of the point on this curve from which the normal is drawn, and $ r $ is the distance from the curve along the normal to it. We assume that $ r> 0 $ in the domain $ \Omega^{(+)}$ and $r <0 $ in the domain $ \Omega^{(-)} $, and note that for $ y = h (l, t) $ we have $ r = 0 $ and $x=l$, and that the derivatives of the functions $ h (x, t) $ in expression \eqref{sincosalpha} are also taken for $ x = l. $

The asymptotic solution to problem  \eqref{mainproblem} will be constructed in the following form:
\begin{equation} \label{asymptoticsolution}
U(x, y, t, \mu)=
\begin{cases}
U^{(-)} (x, y, t, \mu) , \quad (x, y, t) \in \mathbb{R}\times \bar{\Omega}^{(-)} \times \mathcal{T} ,\\
U^{(+)}(x, y, t, \mu) , \quad (x, y,t )\in \mathbb{R}\times\bar{\Omega}^{(+)} \times \mathcal{T},
\end{cases}
\end{equation}
as the sums of two terms
\begin{equation} \label{asymptoticapproximation2terms}
U^{(\mp)}=\bar{u}^{(\mp)}(x, y, \mu)+Q^{(\mp)}(\xi, l, h(l, t), t, \mu),
\end{equation}
where $\bar{u}^{(\mp)}(x, y, t, \mu) $ is the outer functions of the asymptotic representation, and $Q^{(\mp)}(\xi,$ $ l, h(l, t), t, \mu)$ represents the inner functions describing the transition layer. Each term in \eqref{asymptoticsolution} will be represented as an expansion in powers of the small parameter $\mu$:
\begin{align}
\bar{u}^{(\mp)}(x, y, \mu)&=\bar{u}_{0}^{(\mp)}(x, y)+\mu \bar{u}_{1}^{(\mp)}(x, y)+\ldots, \label{expansionregularfunctions} \\
Q^{(\mp)}(\xi, l, h(l, t), t, \mu)&=Q_{0}^{(\mp)}(\xi, l, h(l, t), t)+\mu Q_{1}^{(\mp)}(\xi, l, h(l, t), t)+\ldots. \label{expansiontransitionfunctions}
\end{align}

We assume that $y=h(x, t)$ -- is the curve on which the solution $u(x, y, t, \mu)$ to problem \eqref{mainproblem} at each time instant takes on a value equal to the half-sum of the functions  $\bar{u}^{(-)} (x, y)$ and $\bar{u}^{(+)}(x, y)$:
\begin{equation} \label{halfsum}
\displaystyle u(x, h, t, \mu)=\phi(x, h, \mu):=\frac{1}{2}\left(\bar{u}^{(-)}(x, h,\mu)+\bar{u}^{(+)}(x, h,\mu) \right).
\end{equation}

For the zero-order approximation, it takes the form:
\begin{equation} \label{halfsum0order}
\displaystyle \phi_0(x, h):=\frac{1}{2}\left(\varphi^{(-)}(x, h)+\varphi^{(+)}(x, h) \right).
\end{equation}

The curve $ y = h (x, t) $ will also be sought in the form of an expansion in powers of a small parameter:
\begin{equation} \label{curveexpansion}
 h(x, t)=h_{0}(x, t)+\mu h_{1}(x, t)+\mu^{2}h_{2}(x, t)+\ldots.
\end{equation}

The functions $U^{(-)} (x, y, t, \mu)$ and $U^{(+)}(x, y, t, \mu)$ and their derivatives along the normal to the curve $y=h(x, t)$ are matched on the curve $h(x, t)$ at each moment in time $t$:
\begin{align} 
U^{(-)} (x, h(x, t), t, \mu)=U^{(+)}(x, h(x, t), t, \mu)=\phi(x, h(x, t), t) , \label{sewingcond1} \\
\displaystyle \frac{\partial U^{(-)}}{\partial n}(x, h(x, t), t, \mu)=\frac{\partial U^{(+)}}{\partial n}(x, h(x, t), t, \mu), \label{sewingcond2}
\end{align}
where the function $\phi(x, h(x, t), t)$ is defined in \eqref{halfsum}.

We rewrite the differential operators in equation \eqref{mainproblem} in the variables $r, l, t$:
\begin{align} \label{nabla}
\displaystyle \nabla=\left\lbrace -\frac{h_{l}}{\sqrt{1+h_{l}^{2}}}\frac{\partial}{\partial r}-\frac{\sqrt{1+h_{l}^{2}}}{rh_{ll}-(1+h_{l}^{2})^{\frac{3}{2}}}\frac{\partial}{\partial l};  \frac{1}{\sqrt{1+h_{l}^{2}}}\frac{\partial}{\partial r}-\frac{h_{l}\sqrt{1+h_{l}^{2}}}{rh_{ll}-(1+h_{l}^{2})^{\frac{3}{2}}}\frac{\partial}{\partial l} \right\rbrace ;
\end{align}

\begin{multline} \label{laplasian}
\Delta=\frac{\partial^{2}}{\partial r^{2}}+\frac{(1+h_{l}^{2})^{2}}{(rh_{ll}-(1+h_{l}^{2})^{\frac{3}{2}})^{2}}\frac{\partial^{2}}{\partial l^{2}}+\frac{h_{ll}}{rh_{ll}-(1+h_{l}^{2})^{\frac{3}{2}}}\frac{\partial}{\partial r}\\
+\displaystyle \frac{1+h_{l}^{2}}{(rh_{ll}-(1+h_{l}^{2})^{\frac{3}{2}})^{3}}(2rh_{l}h_{ll}^{2}+h_{l}h_{ll}(1+h_{l}^{2})^{\frac{3}{2}}-rh_{lll}(1+h_{l}^{2}))\frac{\partial}{\partial l}.
\end{multline}

Taking into account that
\begin{equation} \label{xiequation}
\displaystyle \xi=\frac{1}{\mu } \left(y-h(l,t) \right) \sqrt{1+h_{l}^{2}},
\end{equation}
and deriving $l(t) $ from \eqref{localcoord} we rewrite the operator $\partial / \partial t$ in the variables $\xi,l,t$:
\begin{align} \label{operatordt}
\frac{\partial}{\partial t}=-\frac{1}{\mu}\left(h_t \sqrt{1+h_{l}^{2}}-\frac{rh_{l} h_{lt}}{1+h_{l}^{2}} \right) \frac{\partial }{\partial \xi}  -\left(h_t h_l - \frac{r h_{lt} }{\sqrt{1+h_{l}^{2}}} \right)\frac{\partial}{\partial l}+\frac{\partial}{\partial t}. 
\end{align}

We substitute  \eqref{asymptoticapproximation2terms}  into \eqref{mainproblem}, and then subtract the part which describes the outer regions from the equation and replace the variables with $ \xi, l, t $ using equations \eqref{nabla}–\eqref{operatordt}. Following this, we obtain equations for the transition-layer functions  $Q^{(\mp)}(\xi, l, h(l, t), t, \mu)$:
\begin{multline} \label{transitionlayerequation}
\displaystyle  \frac{h_{t}(1+h_{l}^{2})-(h_{l}-1)\left(\bar{u}^{(\mp)}(\xi, l, t, \mu) + Q^{(\mp)}\right)}{\mu \sqrt{1+h_{l}^{2}}}\frac{\partial Q^{(\mp)}}{\partial\xi}
+\frac{1}{\mu} \frac{\partial^{2}Q^{(\mp)}}{\partial\xi^{2}}-\frac{\partial Q}{\partial t}+\left( \frac{\xi h_{l} h_{lt}}{1+h_{l}^{2}} -\frac{h_{ll}}{(1+h_{l}^{2})^{\frac{3}{2}}} \right) \frac{\partial Q^{(\mp)}}{\partial\xi}\\
-\displaystyle \frac{\left(\bar{u}^{(\mp)}(\xi, l, t, \mu)+Q^{(\mp)}\right)(1-h_{l})}{1+h_{l}^{2}} \frac{\partial Q^{(\mp)}}{\partial l}
 +\displaystyle \sum_{i=1}\mu^{i}\mathcal{L}_{i}[Q^{(\mp)}]
 \\= -Q^{(\mp)} \left( \frac{1}{\mu}\frac{1-h_l}{\sqrt{1+h_{l}^{2}}} \frac{\partial \bar{u}^{(\mp)}(\xi, l, t, \mu)}{\partial \xi} +  \frac{h_{l}+1}{1+h_{l}^{2}}\frac{\partial \bar{u}^{(\mp)}(\xi, l, t, \mu)}{\partial l} \right),
\end{multline}
where $ \mathcal{L}_{i} $ represents differential operators of the first or second order in the variables $ \xi $ and $ l $, and the following notations are used:
\begin{align} \label{regularfunctionschangevariable}
\bar{u}^{(\mp)}(\xi, l, t, \mu) := \bar{u}^{(\mp)}(l-\mu\xi\sin\alpha, h(l, t)+\mu\xi\cos\alpha).
\end{align}

\subsection{Zero-order approximation terms}

We find the outer functions of zero order by substituting expansions \eqref{expansionregularfunctions} into the stationary equation
\begin{equation} \label{ArbRegFuncEq}
    \displaystyle \mu  \Delta \bar{u} =-\bar{u} \left( k \frac{\partial \bar{u}}{\partial x}+ \frac{\partial \bar{u}}{\partial y} \right)+f(x, y).
\end{equation}
Expanding the functions into the series in powers of a small parameter and equating the coefficients of $ \mu^{0} $, we obtain the following reduced equations:
\begin{align}\label{zeroorderregularequation1}
\begin{cases}
\displaystyle \bar{u}_{0}^{(-)} \left(k\frac{\partial\bar{u}_{0}^{(-)}}{\partial x}+\frac{\partial\bar{u}_{0}^{(-)}}{\partial y}\right)=f(x, y),\\
\bar{u}_{0}^{(-)}(x, -a)=u^{-a}(x), \quad \bar{u}_{0}^{(-)}(x, y)=\bar{u}_{0}^{(-)}(x+L, y);
\end{cases}
\end{align}
\begin{align}\label{zeroorderregularequation2}
\begin{cases}
\displaystyle \bar{u}_{0}^{(+)}\left( k\frac{\partial\bar{u}_{0}^{(+)}}{\partial x}+\frac{\partial\bar{u}_{0}^{(+)}}{\partial y}\right)=f( x, y),\\
\bar{u}_{0}^{(+)}(x, a)=u^{a}(x), \quad \bar{u}_{0}^{(+)}(x, y)=\bar{u}_{0}^{(+)}(x+L, y).
\end{cases}
\end{align}

According to Assumption \ref{A1}, equations \eqref{zeroorderregularequation1} and \eqref{zeroorderregularequation2} have the following solutions:
\begin{align}
\label{u0regu}
\bar{u}_{0}(x,y)= \begin{cases}
\bar{u}_{0}^{(-)}(x,y)=\varphi^{(-)}(x,y), \quad (x,y)\in \mathbb{R} \times \bar{\Omega}^{(-)},\\
\bar{u}_{0}^{(+)}(x,y)=\varphi^{(+)}(x,y), \quad (x,y)\in  \mathbb{R} \times \bar{\Omega}^{(+)}.
\end{cases}
\end{align}

We substitute the series \eqref{expansionregularfunctions}, \eqref{expansiontransitionfunctions}, and \eqref{curveexpansion} into equations \eqref{sewingcond1} and \eqref{transitionlayerequation}, expand all the terms of \eqref{transitionlayerequation} in series in powers of $\mu$, equate the coefficients for $\mu^{-1}$ in equation \eqref{transitionlayerequation} and those for $\mu^{0}$ in equation \eqref{sewingcond1}, and taking into account the additional condition for the decay of the transition functions at infinity, we obtain the problems for transition-layer functions of zero order $Q_{0}^{(\mp)}(\xi, l, h_0(l, t), t)$:
\begin{align} \label{transitionalfunczeroord}
\begin{cases}
\displaystyle \frac{\partial^{2}Q_{0}^{(\mp)}}{\partial\xi^{2}}
+\displaystyle \frac{{h_0}_{t}\left(1+{h_0}_{l}^{2}\right)+\left(1-k {h_0}_{l}\right)\left( \varphi^{(\mp)}(l, h_0)+Q_{0}^{(\mp)} \right)}{\sqrt{1+{h_0}_{l}^{2}}}\frac{\partial Q_{0}^{(\mp)}}{\partial\xi}=0,\\
\varphi^{(\mp)}(l, h_0)+Q_{0}^{(\mp)}(0, l, h_0, t)=\phi_0(l, h_0),\\
Q_{0}^{(\mp)}(\mp\infty, l, h_0, t)=0,
\end{cases}
\end{align}
where $\phi_0(l, h_0)$ is obtained from \eqref{halfsum0order}.

\begin{remark}
To determine $Q_{0}^{(-)}$ in the problem \eqref{transitionalfunczeroord}, we consider values $\xi\leq 0$, while to determine $Q_{0}^{(+)}$, we consider values $ \xi\geq 0$.
\end{remark}

To study the zero approximation of $h(l,t)$, i.e. $h_0(l,t)$, we introduce the auxiliary function
\begin{equation} \label{auxilaryfunction}
\tilde{u}(\xi, h_0(l, t))=
\begin{cases}
\varphi^{(-)} (l, h_0(l, t))+Q_{0}^{(-)} (\xi, l, h_0(l, t), t) , \quad \xi\leq 0,\\
\varphi^{(+)}(l, h_0(l, t))+Q_{0}^{(+)}(\xi, l, h_0(l, t), t) , \quad \xi\geq 0.
\end{cases}
\end{equation}
Each of the separate equations in \eqref{transitionalfunczeroord}, written in the above notations, takes the following form:
\begin{equation} \label{replacement}
\displaystyle \frac{\partial^{2}\tilde{u}}{\partial\xi^{2}}+\frac{{h_0}_{t}(l, t)\left(1+{h_0}_{l}^{2}(l,t)\right)+\tilde{u}\left(1-k {h_0}_{l}(l, t)\right)}{\sqrt{1+{h_0}_{l}^{2}(l,t)}}\frac{\partial\tilde{u}}{\partial\xi}=0.
\end{equation}

 Let $\displaystyle \frac{\partial \tilde{u} }{\partial\xi} = g(\tilde{u}), \frac{\partial^2 \tilde{u} }{\partial\xi^{2}} = \frac{\partial g(\tilde{u}) }{\partial \tilde{u} } g(\tilde{u}) $; then, \eqref{replacement} is transformed into
$$\displaystyle \frac{\partial g(\tilde{u}) }{\partial \tilde{u} } = \frac{-{h_0}_{t}(l, t)\left(1+{h_0}_{l}^{2}(l,t)\right)+\tilde{u}\left(k{h_0}_{l}(l, t)-1\right)}{\sqrt{1+{h_0}_{l}^{2}(l,t)}},$$
from which we can deduce that
\begin{align} \label{derivativetildeu}
\frac{\partial\tilde{u}}{\partial\xi}= \begin{cases}
\displaystyle \Phi^{(-)} (\xi, h_0)=\int_{\varphi^{(-)}}^{\tilde{u}} \frac{-{h_0}_{t}\left(1+{h_0}_{l}^{2}\right)+u\left(k{h_0}_{l}-1\right)}{\sqrt{1+{h_0}_{l}^{2}}} du , \quad \xi\leq 0,\\
\displaystyle \Phi^{(+)}(\xi, h_0)= \int_{\varphi^{(+)}}^{\tilde{u}} \frac{-{h_0}_{t}\left(1+{h_0}_{l}^{2}\right)+u\left(k{h_0}_{l}-1\right)}{\sqrt{1+{h_0}_{l}^{2}}} du , \quad \xi\geq 0.
\end{cases}
\end{align}

%We rewrite matching conditions \eqref{sewingcond2}, writing the derivative along the normal to the curve $h(x, t)$; in variables $x, y, t$, we have
Let us write the derivative along the normal to the curve $h(x, t)$ with respect to the variables $x, y, t$:
\begin{equation} \label{NormalToCurveEq}
\displaystyle \frac{\partial}{\partial n}=(n, \nabla)=\frac{\partial}{\partial r}=-\sin\alpha\frac{\partial}{\partial x}+\cos\alpha\frac{\partial}{\partial y},  \end{equation}
where $\sin\alpha$ and $\cos\alpha$ are determined by expression \eqref{sincosalpha}.

In the variables $\xi, l, t$ this derivative has the following form:
\begin{equation} \label{derivativetonormal}
    \frac{\partial}{\partial n}=\frac{1}{\mu} \frac{\partial}{\partial \xi}.
\end{equation}

Taking into account expansions \eqref{asymptoticapproximation2terms}--\eqref{curveexpansion} and \eqref{derivativetonormal}, and equating the coefficients at $\mu^{-1}$ in \eqref{sewingcond2}, we obtain (assuming that at $\xi=0, x=l$ and $y=h_0(x,t)$):
\begin{align}\label{sewindcondexpanded0}
\Phi^{(-)} (0, h_0(x, t))-\Phi^{(+)} (0, h_0(x, t))=0.
\end{align}

From \eqref{derivativetildeu} and \eqref{sewindcondexpanded0}, we obtain
\begin{equation*}
\int_{\varphi^{(-)} (x, h_0(x, t))}^{\varphi^{(+)}(x, h_0(x, t))} \frac{-{h_0}_{t}(x, t)\left(1+{h_0}_{x}^{2}(x,t)\right)+u\left(k{h_0}_{x}(x, t)-1\right)}{\sqrt{1+{h_0}_{x}^{2}(x,t)}} du =0.
\end{equation*}

From this, we can find the equation that determines the zero approximation of the curve $h_0(x,t)$:
\begin{equation} \label{h0mainequation}
{h_0}_{t}\left(1+{h_0}_{x}^{2}\right)=\frac{1}{2} \left(k{h_0}_{x}-1\right) \left(\varphi^{(+)}(x, h_0)+\varphi^{(-)}(x, h_0) \right),
\end{equation}
with conditions ${h_0}(x, t)={h_0}(x+L, t), {h_0}(x, 0)=h_{0}^{*} \in (-a;a) $.

We solve equation \eqref{derivativetildeu}, and, with \eqref{h0mainequation}, the functions $Q_{0}^{(\mp)}(\xi, l, h_0, t)$ can be represented in the explicit form  in which $ h_0 (x, t) $ is a parameter:
\begin{align} \label{Q0equation}
\displaystyle Q_{0}^{(\mp)}(\displaystyle \xi, l, h_0, t)=\frac{2P^{(\mp)}(x,h_0)}{\exp \left(-\xi \frac{  P^{(\mp)}(x,h_0) \left(1-k{h_0}_{x} \right)}{\sqrt{1+{h_0}_{x}^{2}}} \right)+1},
\end{align}
where
\begin{align*}
 P^{(-)}(x,h_0)=\frac{1}{2}\left(\varphi^{(+)}(x,h_0)-\varphi^{(-)}(x,h_0) \right), \quad
 P^{(+)}(x,h_0)=\frac{1}{2}\left(\varphi^{(-)}(x,h_0)-\varphi^{(+)}(x,h_0)\right).
\end{align*}

Consequently, the transition-layer functions $Q_{0}^{(\mp)}(\displaystyle \xi, l, h_0(l, t), t)$ are exponentially decreasing with $\xi \rightarrow \mp \infty$ and have the exponential estimates \cite{Vasileva1998ContrastSI,Butuzov1997ASYMPTOTICTO}:
\begin{equation}\label{equat22}
\undersl{C} e^{\undersl{\kappa}\xi} \leq |Q_{0}^{(-)}( \xi, l, h_0(l, t), t)|\leq \oversl{C}e^{\oversl{\kappa}\xi} , \quad \xi\leq 0, \quad t\in \bar{\mathcal{T}},
 \end{equation}
\begin{equation}\label{equat23}
\undersl{C} e^{-\undersl{\kappa}\xi} \leq |Q_{0}^{(+)}( \xi, l, h_0(l, t), t)|\leq \oversl{C}e^{-\oversl{\kappa}\xi}, \quad \xi\geq 0, \quad t\in \bar{\mathcal{T}},
 \end{equation}
where $\undersl{C}, \oversl{C}$ and $\undersl{\kappa}, \oversl{\kappa}$ -- are four positive constants independent of $\xi,l,t$, and, precisely, $$\undersl{C} := \frac{1}{2}  \inf_{t \in [0,T]} \varphi^{(+)}(x,h_{0}(x,t))-  \varphi^{(-)}(x,h_{0}(x,t)), $$
$$ \oversl{C} := \frac{1}{2}  \sup_{t \in [0,T]} \varphi^{(+)}(x,h_{0}(x,t))-  \varphi^{(-)}(x,h_{0}(x,t)).$$

From the boundary conditions of \eqref{transitionalfunczeroord}, and from Assumption  \ref{A1}, we deduce that
$$|Q_{0}^{( \mp)}( 0, x, h_0(x, t), t)|=|\phi_0((x, h_0(x, t)))-\varphi^{(\mp)}((x, h_0(x, t)))| > \frac{1}{2} \left( u^a(x)-u^{-a}(x) \right) > \mu^2 $$
and $\lvert Q_{0}^{(\mp)}( \xi, l, h_0(l, t), t) \rvert \rightarrow 0$  for $\xi \rightarrow \mp \infty$. Since $\mu >0 $ is a constant, $\xi=\frac{1}{\mu}  \frac{y-h(l, t)}{\cos{\alpha}}$  and  $|Q_{0}^{(\mp)}( \xi, l,$ $ h_0(l, t), t)|$ is a decreasing function, and there exists $y=H^{(\mp)}(l, t)$   starting from which we have  $|Q_{0}^{(\mp)}(\xi, l, H^{(\mp)}, t)| \leq \mu^2 $ for every $l,t$, i.e.
\begin{equation}  \label{Qlr-mu}
|Q_{0}^{(\mp)}(\xi( H^{(\mp)}), l, H^{(\mp)}, t)| = \mu^2.
\end{equation}

We denote the width of the transition layer with $ \Delta h = H^{(+)}(l, t) -H^{(-)}(l, t)$, and from equations \eqref{equat22}–\eqref{equat23}, for $y=H^{(\mp)}(l, t)$, we obtain
\begin{equation}
\displaystyle \undersl{C} e^{-\frac{\underline{\kappa}}{2\mu  \cos{\alpha} } \Delta h} \leq \mu^2 \leq \oversl{C}e^{-\frac{\bar{\kappa}}{2\mu  \cos{\alpha}} \Delta h},
\end{equation}
from which $ \Delta h$ can be estimated as
\begin{equation} \label{Deltah}
\frac{2\mu  \cos{\alpha}}{\undersl{\kappa}} \ln \frac{\undersl{C}}{\mu^2} \leq \Delta h \leq \frac{2\mu  \cos{\alpha}}{\oversl{\kappa}} \ln \frac{\oversl{C}}{\mu^2}, \text{~i.e.~} \Delta h \sim \mu |\ln \mu|.
\end{equation}

\subsection{First-order approximation terms}

First-order-approximation outer functions $\bar{u}_{1}^{(\mp)}(x,y)$ are defined from \eqref{ArbRegFuncEq} as solutions to problems
\begin{equation} \label{firstorderregularequation}
\begin{split}
&k \frac{\partial\bar{u}_{1}^{(\mp)}}{\partial x} + \frac{\partial\bar{u}_{1}^{(\mp)}}{\partial y}+\bar{u}_{1}^{(\mp)} P(x,y)=W(x,y),\\ &\bar{u}_{1}^{(-)}(x, -a)=0,  \quad \bar{u}_{1}^{(-)}(x, y)=\bar{u}_{1}^{(-)}(x+L, y), \\ &\bar{u}_{1}^{(+)}(x, a)=0, \quad \bar{u}_{1}^{(+)}(x, y)=\bar{u}_{1}^{(+)}(x+L, y),
\end{split}
\end{equation}
where
$$P(x,y)= \frac{1}{\varphi^{(\mp)}} \left( k \frac{\partial\varphi^{(\mp)}}{\partial x} + \frac{\partial\varphi^{(\mp)}}{\partial y}\right), \  W(x,y)= -\frac{1}{\varphi^{(\mp)}} \left( \frac{\partial^{2}\varphi^{(\mp)}}{\partial x^{2}}+\frac{\partial^{2}\varphi^{(\mp)}}{\partial y^{2}} \right). $$

From \eqref{firstorderregularequation}, we obtain  $\bar{u}_{1}^{(\mp)}(x, y)$ in its explicit form:
\begin{align}\label{firstoredregularfunctions}
\bar{u}_{1}^{(\mp)}=\displaystyle \exp \left(\int _1^x-\frac{P\left(s,\frac{s-x+k y}{k}\right)}{k}ds\right)  \int _{ \mp ak +x-k y}^x\frac{\exp \left(\int _1^{\zeta }\frac{P\left(s,\frac{s-x+k y}{k}\right)}{k}ds\right) W\left(\zeta ,\frac{-x+k y+\zeta }{k}\right)}{k}d\zeta.
\end{align}

The equations for the first-order transition-layer functions $Q_{1}^{(\mp)}(\xi, l, h_0(l, t), t)$ are obtained by equating terms at $\mu^{0}$ in \eqref{transitionlayerequation}:
\begin{multline}\label{transitionlayerfirstord}
\displaystyle \frac{\partial^{2}Q_{1}^{(\mp)}}{\partial\xi^{2}}+  \frac{\partial}{\partial \xi} \left(Q_{1}^{(\mp)} \frac{{h_0}_{t}  ({h_0}_{l}^{2} +1  ) + \tilde{u}(1-k{h_0}_{l} ) }{ \sqrt{1+{h_0}_{l}^{2} }}   \right)= r_{1}^{(\mp)}(\xi,l,t) \frac{\partial h_{1} }{\partial l}\\
+r_{2}^{(\mp)}(\xi,l,t) h_{1} +r_{3}^{(\mp)}(\xi,l,t)\frac{\partial h_{1} }{\partial t}+r_{4}^{(\mp)}(\xi,l,t):=f_{1}^{(\mp)}(\xi,l,t),
\end{multline}
where $r_{1}^{(\mp)}(\xi,l,t), r_{2}^{(\mp)}(\xi,l,t), r_{3}^{(\mp)}(\xi,l,t), r_{4}^{(\mp)}(\xi,l,t)$ are known functions, and, particularly, $\displaystyle r_{3}^{(\mp)}(\xi,l,t)$ $=$ $\frac{\partial Q_{0}^{(\mp)} }{\partial \xi}  \sqrt{1+{h_0}_{l}^{2}}$. The derivatives of the function $ h_1 (l, t) $ are taken for $ x = l $. From \eqref{sewingcond1}, taking into account the matching conditions of zero order in \eqref{transitionalfunczeroord} and \eqref{halfsum}, we obtain the boundary conditions
\begin{align}\label{transitionlayerfirstordbound1}
\begin{split}
Q_{1}^{(\mp)}(0, x, h_0, t)=\frac{\bar{u}_{1}^{(\pm)}(x, h_0)-\bar{u}_{1}^{(\mp)}(x, h_0)}{2}:= p_{1}^{(\mp)} (x, h_0).
 \end{split}
\end{align}
We also add conditions at infinity:
\begin{equation}\label{transitionlayerfirstordbound2}
Q_{1}^{(\mp)} (\mp \infty, l, h_0(l, t), t)=0.
\end{equation}
The solutions to problems \eqref{transitionlayerfirstord}–\eqref{transitionlayerfirstordbound2} can be written explicitly:
\begin{align} \label{Q1function}
Q_{1}^{(\mp)}(\xi, l, h_0, t) 
=J^{(\mp)}(\xi,h_0) \left( p_{1}^{(\mp)}(x, h_0)  \vphantom{\int_{0}^{0}}  +\int_{0}^{\xi}\frac{ds}{J^{(\mp)}(s,h_0)}\int_{\mp\infty}^{s}f_{1}^{(\mp)}(\eta, l, t)d\eta ds \right),
\end{align}
where $ \displaystyle J^{(\mp)}(\xi,h_0)= \left( \Phi^{(\mp)}(0,h_0) \right)^{-1} \Phi^{(\mp)}(\xi,h_0). $

From \eqref{Q1function}, we find:
\begin{align} \label{Q1functionDerivative}
\frac{\partial Q_{1}^{( \mp)}}{\partial \xi} (0, x, h_0, t) 
=p_{1}^{(\mp)}(x, h_0) \left( -{h_0}_{t}\sqrt{1+{h_0}_{x}^{2}}  \vphantom{\frac{\phi_0(x, h_0)\left(k{h_0}_{x}(x, t)-1\right)}{\sqrt{1+{h_0}_{x}^{2}(x,t)}}} \right. \left. +\frac{\phi_0(x, h_0)\left(k{h_0}_{x}-1\right)}{\sqrt{1+{h_0}_{x}^{2}}} \right)  - \int_{0}^{\mp \infty} f_{1}^{( \mp)} (\eta,x,t) d\eta.
\end{align}

It can be clearly shown that $Q_{1}^{( \mp)} (\xi, l, h_0(x, t), t) $ satisfies exponential estimates \eqref{equat22} and \eqref{equat23}. From the first-order $C^1$-matching condition \eqref{sewingcond2}, and taking into account expansions \eqref{sincosalpha}, \eqref{NormalToCurveEq}, and \eqref{derivativetonormal}, we obtain
\begin{multline} \label{matchingfirstord}
 h_1 \left(\frac{\partial^2 Q_{0}^{(-)}}{\partial \xi \partial h_0}-\frac{\partial^2 Q_{0}^{(+)}}{\partial \xi \partial h_0} \right)+\frac{\partial Q_{1}^{(-)}}{\partial \xi }-\frac{\partial Q_{1}^{(+)}}{\partial \xi } \\
 +\frac{1}{\sqrt{1+{h_0}_{x}^{2}}} \left( \frac{\partial \varphi^{(-)}}{\partial h_0}-\frac{\partial \varphi^{(+)}}{\partial  h_0} \right)
 +\frac{{h_0}_{x}}{\sqrt{1+{h_0}_{x}^{2}}} \left( \frac{\partial \varphi^{(-)}}{\partial x}-\frac{\partial \varphi^{(+)}}{\partial x} \right)=0.
\end{multline}
From \eqref{transitionlayerfirstordbound1}, \eqref{Q1functionDerivative}, and \eqref{matchingfirstord}, we obtain the equation determining $h_1(x,t)$:
\begin{align}\label{eqforh1}
\begin{split}
  & {h_1}_{t}  \left( \varphi^{(+)}-\varphi^{(-)} \right)\sqrt{1+{h_0}_{x}^{2}} ={h_1}_{x} V_1(x,t) + h_1 V_2(x,t) +V_3(x,t), \\
  & h_{1}(x, t)=h_{1}(x+L, t),\ h_1(x,0)=0,
\end{split}
\end{align}
where $ V_1(x,t),V_2(x,t),V_3(x,t)$ depend on known functions.

Problem \eqref{eqforh1} is solvable because the coefficient of the term ${h_1}_{t}$ is positive.

In a similar way to \eqref{firstorderregularequation}–\eqref{eqforh1}, we can obtain approximation terms for the asymptotic solution up to an arbitrary order.
We write the approximation terms of $h(x, t)$ in \eqref{curveexpansion} up to order $n$:
\begin{equation} \label{hExpansion}
\hat{h}_{n}(x,t)=\sum_{i=0}^{n} \mu^{i} h_{i}(x,t), \quad x \in \mathbb{R}, \ t \in \bar{\mathcal{T}}.
\end{equation}

\section{Main results} \label{statementresults}

The notations and abbreviations frequently used in this section are listed in Table \ref{NotationTable} below.

\begin{longtable}{l l l}
\hline \multicolumn{1}{c}{\textbf{Notation}} & \multicolumn{1}{c}{\textbf{Description }} & \multicolumn{1}{c}{\textbf{Reference}} \\ \hline
\endfirsthead
\hline
\hspace{-3.5mm} \begin{tabular}{l}  $\mu$   \end{tabular} & \hspace{-3.5mm} \begin{tabular}{l}  Small parameter, $ 0<\mu \ll 1 $       \end{tabular}   &  Eq.\eqref{mainproblem}  \\
 \hline
\renewcommand{\arraystretch}{1.2} \hspace{-4.5mm}  \begin{tabular}{l} $\Omega$, $\bar{\mathcal{T}}$ \end{tabular} & \hspace{-3.5mm} \begin{tabular}{l} Domains: $\Omega=(-a,a)$, $\bar{\mathcal{T}}=[0, T]$ \end{tabular} & Eq.\eqref{mainproblem} \\
 \hline
\hspace{-3.5mm} \begin{tabular}{l}    $u^{(\mp)}$ \end{tabular} & \hspace{-3.5mm} \begin{tabular}{l}  Left and right boundary conditions \end{tabular}  & Eq.\eqref{mainproblem} \\
   \hline
\hspace{-3.5mm} \begin{tabular}{l}  $h(x,t)$ \end{tabular} & \hspace{-3.5mm} \begin{tabular}{l}  The curve that is located in the middle \\ of the transition layer  \end{tabular}  & Eq.\eqref{curveexpansion} \\
 \hline
\hspace{-3.5mm} \begin{tabular}{l}  $\hat{h}_n(x,t)$ \end{tabular} & \hspace{-3.5mm} \begin{tabular}{l} The expansion of  $h(x,t)$ in a series \\ up to order $n$   \end{tabular}  & Eq.\eqref{hExpansion} \\
 \hline
 \hspace{-3.5mm} \begin{tabular}{l}  $\displaystyle \frac{\partial }{\partial n}$ \end{tabular} & \hspace{-3.5mm} \begin{tabular}{l}  Derivative along the normal \\ to the curve $h(x,t)$  \end{tabular}  & Eq.\eqref{NormalToCurveEq} \\
 \hline
\hspace{-3.5mm} \begin{tabular}{l} Superscript  $ ^{(\mp)}$ \end{tabular} & \hspace{-3.5mm} \begin{tabular}{l} Describe functions on the left and right, \\
respectively, relative to the curve $h(x,t)$
\end{tabular} & Eq.\eqref{asymptoticsolution}   \\
   \hline
\hspace{-3.5mm} \begin{tabular}{l} Subscript   $ _{0,1...}$ \end{tabular} & \hspace{-3.5mm} \begin{tabular}{l}  The order of approximation of the \\ asymptotic solution \end{tabular} & Eq.\eqref{expansionregularfunctions}-\eqref{curveexpansion} \\
 \hline
\hspace{-4.5mm} \renewcommand{\arraystretch}{1.2}  \begin{tabular}{l}  $\bar{\Omega}^{(-)}$, $\bar{\Omega}^{(+)}$ \end{tabular} & \hspace{-3.5mm} \begin{tabular}{l}  $[-a,h(x,t)]$ and  $[h(x,t),a]$ respectively     \end{tabular}   & Eq.\eqref{asymptoticsolution}  \\
 \hline
\hspace{-3.5mm} \begin{tabular}{l} $\bar{u}^{(\mp)}$ \end{tabular} & \hspace{-3.5mm} \begin{tabular}{l} Outer functions describing the   \\  solution far from the curve $h(x,t)$ \end{tabular} & Eq.\eqref{asymptoticapproximation2terms}  \\
 \hline
\hspace{-3.5mm} \begin{tabular}{l} $\varphi^{(\mp)}$ \end{tabular} & \hspace{-3.5mm} \begin{tabular}{l} Zero-order approximation \\ of outer functions \end{tabular} & Eq.\eqref{eq7},\eqref{u0regu}  \\
  \hline
\hspace{-3.5mm} \begin{tabular}{l} $Q^{(\mp)}$ \end{tabular} &  \hspace{-3.5mm} \begin{tabular}{l} Transition-layer functions describing \\ the solution  near the curve $h(x,t)$ \end{tabular} & Eq.\eqref{expansiontransitionfunctions} \\
 \hline
\hspace{-3.5mm}  \begin{tabular}{l}  $ \Delta h$  \end{tabular} & \hspace{-3.5mm} \begin{tabular}{l}   The width of the transition layer, \\ $ \Delta h \sim \mu |\ln \mu|$  \end{tabular}   &  Eq.\eqref{Deltah} \\
 \hline
\hspace{-3.5mm} \begin{tabular}{l} $ \xi$   \end{tabular} & \hspace{-3.5mm} \begin{tabular}{l}  Extended variable, $\xi = r/ \mu$       \end{tabular}   & Eq.\eqref{xidefin} \\
 \hline
  \hspace{-4.5mm}  \renewcommand{\arraystretch}{1.2}
  \begin{tabular}{l} $u^\varepsilon (x,y,t)$ \end{tabular}  & \hspace{-3.5mm} \begin{tabular}{l}     The smooth approximate data \end{tabular}   & Eq.\eqref{uAlphaL},\eqref{uAlphaR} \\
 \hline
  \hspace{-4.5mm}  \renewcommand{\arraystretch}{1.2}
  \begin{tabular}{l} $f^*(x,y) $ \end{tabular}  & \hspace{-3.5mm} \begin{tabular}{l}      The exact source function   \end{tabular}   & Eq.\eqref{f0Ineq} \\
 \hline
 \hspace{-4.5mm}  \renewcommand{\arraystretch}{1.2}
  \begin{tabular}{l} $f^\delta(x,y,t) $ \end{tabular}  & \hspace{-3.5mm} \begin{tabular}{l}      The regularized approximate \\ source function \end{tabular}   & Eq.\eqref{fdelta} \\
 \hline
\caption{Notations and references to their definitions.}
\label{NotationTable}
\end{longtable}

We are now able to determine the main result for the forward problem (the meanings and descriptions of some notations can be found in Table \ref{NotationTable} and Section \ref{constructionOfAsymptotic}).

\begin{theorem} \label{MainThm}
Suppose that functions $f(x,y), u^{-a}(x), u^{a}(x), u_{init}(x,y)$ are sufficiently smooth and $x$-periodic, and $\mu\ll 1$. Then, under Assumptions \ref{A1}–\ref{A4}, the boundary-value problem \eqref{mainproblem} has a unique smooth solution with an inner transition layer. In addition, the $n$-order asymptotic solution $U_{n}(x,y,t,\mu)$ has the following representation ($\xi_{n} = \left(y-\hat{h}_n \right) \sqrt{1+(\hat{h}_n)_{x}^{2}}/ \mu$):
\begin{align} \label{asymptoticnorder}
 U_{n}= \begin{cases}
 \displaystyle U_{n}^{(-)}=\sum_{i=0}^{n} \mu^{i} \left(\bar{u}_{i}^{(-)}\left(x,y\right)+Q_{i}^{(-)}\left(\xi_{n},l,\hat{h}_{n},t \right) \right), \ (x,y,t) \in \mathbb{R} \times \bar{\Omega}^{(-)} \times \bar{\mathcal{T}} , \\
 \displaystyle  U_{n}^{(+)}=\sum_{i=0}^{n} \mu^{i} \left(\bar{u}_{i}^{(+)}\left(x,y\right)+Q_{i}^{(+)}\left(\xi_{n},l,\hat{h}_{n},t\right)\right), \ (x,y,t) \in \mathbb{R} \times \bar{\Omega}^{(+)}\times \bar{\mathcal{T}}.
  \end{cases}
 \end{align}

Moreover, the following asymptotic estimates hold:
\begin{equation} \label{NorderEstim1}
\forall (x,y,t)\in \mathbb{R} \times \bar{\Omega}  \times \bar{\mathcal{T}}:~ |u(x,y,t)-U_{n}(x,y,t,\mu)|=\mathcal{O}(\mu^{n+1}) ,
\end{equation}
\begin{equation} \label{NorderEstim2}
\forall (x,t) \in \mathbb{R} \times \bar{\mathcal{T}}:~ |h(x,t)-\hat{h}_{n}(x,t)|=\mathcal{O}(\mu^{n+1}),
\end{equation}
\begin{equation}\label{NorderEstim3}
\forall (x,y,t) \in \mathbb{R} \times \bar{ \Omega} \backslash \{ \hat{h}_{n}(x,t) \}  \times \bar{\mathcal{T}}:~ \left| \frac{ \partial u(x,y,t)}{ \partial n}-\frac{ \partial  U_{n}(x,y,t,\mu) }{ \partial n} \right| =\mathcal{O}(\mu^{n}).
\end{equation}
\end{theorem}

\begin{corollary}
\label{Corollary1}
(Zeroth approximation) Under the assumptions of Theorem \ref{MainThm},  the zero-order asymptotic solution $U_{0}$ has the following representation:
 \begin{align} \label{eq001}
U_{0}(x,y,t)=\begin{cases}
\varphi^{(-)}(x,y)+Q_{0}^{(-)}(\xi_{0},l,h_{0}(l,t),t) ,  \ (x,y,t)\in \mathbb{R} \times\bar{\Omega}^{(-)}  \times \bar{\mathcal{T}},\\
\varphi^{(+)}(x,y)+Q_{0}^{(+)}(\xi_{0},l,h_{0}(l,t),t) , \  (x,y,t)\in \mathbb{R} \times\bar{\Omega}^{(+)}  \times \bar{\mathcal{T}},
\end{cases}
\end{align}
where $\xi_{0}=\left(y-h_0 \right) \sqrt{1+(h_0)_{x}^{2}}/ \mu$. Moreover, the following holds:
\begin{equation} \label{eq002}
\forall (x,y,t)\in \mathbb{R} \times \bar{\Omega}  \times \bar{\mathcal{T}}:~ |u(x,y,t)-U_{0}(x,y,t)|=\mathcal{O}(\mu) ,
\end{equation}
\begin{equation} \label{eq003}
\forall t \in \bar{\mathcal{T}}:~ |h(x,t)-h_{0}(x,t)|=\mathcal{O}(\mu) .
\end{equation}

Furthermore, outside the narrow region $\left(h_{0}\left(x,t\right)- \Delta h /2, h_{0}(x,t)+\Delta h /2\right)$ with $\Delta h\sim\mu|\ln\mu|$, there exists a constant $C$ independent of $x, y, t, \mu$ such that the following inequalities hold:
\begin{equation} \label{eq004}
|u(x,y,t)-\varphi^{(-)}(x,y)| \leq C \mu, \qquad  (x,y,t) \in \mathbb{R} \times [-a, h_{0}(x,t)- \Delta h /2]  \times \bar{\mathcal{T}},
\end{equation}
\begin{equation} \label{eq005}
|u(x,y,t)-\varphi^{(+)}(x,y)|  \leq C \mu, \qquad  (x,y,t) \in \mathbb{R} \times [h_{0}(x,t)+\Delta h /2, a]  \times \bar{\mathcal{T}},
\end{equation}

\begin{equation}\label{0orderEstim1}
\left| \frac{\partial u(x,y,t)}{\partial n}-\frac{\partial  \varphi^{(-)}(x,y)}{\partial n} \right|  \leq C \mu,  \quad  (x,y,t) \in \mathbb{R} \times [-a, h_{0}(x,t)- \Delta h /2]   \times \bar{\mathcal{T}},
\end{equation}
\begin{equation}\label{0orderEstim2}
\left| \frac{\partial u(x,y,t)}{\partial n}-\frac{\partial  \varphi^{(+)}(x,y)}{\partial n} \right| \leq C \mu,  \quad  (x,y,t) \in \mathbb{R} \times [h_{0}(x,t)+\Delta h /2, a]  \times \bar{\mathcal{T}}.
\end{equation}
\end{corollary}

Corollary \ref{Corollary1} follows directly from Theorem \ref{MainThm}. Inequalities \eqref{0orderEstim1}–\eqref{0orderEstim2} in Corollary \ref{Corollary1} can be obtained by taking into account the fact that  the transition-layer functions are decreasing functions with respect to the $\xi_0$ and are sufficiently small at the boundaries of the narrow region $(h_0(x,t)-\Delta h/2, h_0(x,t)+\Delta h/2)$, i.e. equation \eqref{Qlr-mu}.

From Corollary \ref{Corollary1}, it follows that the solution can be approximated by outer functions of zero order everywhere, except for a thin transition layer. In view of this, we construct the approximate source function $f$ by using only the outer functions of zero order, which can reduce the computational cost significantly but keep the accuracy of the result in case of a relatively high noise level and small size $\mu$. To do this, suppose we have the deterministic noise model
\begin{equation}
\label{noisyData1}
\max \left( \lvert u(x_i,y_j,t_0)-u^{\delta}_{i,j} \rvert, \left \lvert  u_x(x_i,y_j,t_0)-({u_x})^{\delta}_{i,j} \right \rvert, \left \lvert  u_y(x_i,y_j,t_0)-({u_y})^{\delta}_{i,j} \right \rvert \right) \leq \delta
\end{equation}
between the noisy data $\{u^{\delta}_{i,j}, ({u_x})^{\delta}_{i,j}, ({u_y})^{\delta}_{i,j} \}$ and the corresponding exact quantities $\{u(x_i,$ $y_j,t_0)$,  $u_x(x_i,y_j,t_0)$, $u_y (x_i,y_j,$ $t_0)\}$ at time $t_0$  and at grid points $X_n:= \{x_0 < x_1 < \cdots < x_n \}  $ and $Y_m := \{-a=y_0 < y_1 < \cdots < y_m=a \}$, with maximum mesh sizes in each direction: $d_1:= \max\limits_{ i\in \{0, \cdots,n-1\}  } \{x_{i+1} - x_{i}\}$ and $d_2:= \max\limits_{ j\in \{0, \cdots,m-1\} } \{  y_{j+1} - y_{j}\}$.

For the error estimation, we define the pre-approximate source function $f_0$ 
\begin{equation} \label{f0}
f_0(x,y) =  u (x,y,t_0) \left( k \frac{\partial u (x,y,t_0)}{\partial x} +\frac{\partial u (x,y,t_0)}{\partial y} \right).
\end{equation}

Then, the following assertion holds for $f_0$.
\begin{proposition}
\label{ProAsympErr}
Let $f^* $ be the exact source function, satisfying the original governing equation \eqref{mainproblem}. Under Assumption \ref{A2},
there exists a constant $C_1$ independent of $x,y,\mu$ such that, for any $ p\in(0,+\infty) $,
\begin{equation}
\label{f0Ineq}
\|f^* - f_0 \|_{L^p(\mathbb{R} \times\Omega)} \leq C_1 \mu |\ln \mu | .
\end{equation}
\end{proposition}

According to Corollary \ref{Corollary1}, we can exclude data values $\{u^{\delta}_{i,j}, ({u_x})^{\delta}_{i,j}, ({u_y})^{\delta}_{i,j} \}$  belonging to the transition layer, and  use only nodes from the two regions $[-a, h_0 (x,t)-\Delta h/2] $ and $[ h_0 (x,t)+\Delta h/2,a] $ with node indices $j\in \mathcal{J}$, where $\mathcal{J} := \{ 0, \cdots, m^{(-)},$ $m^{(+)}, \cdots, m \}$. We also introduce notations for the left region, $\bar{\Omega}^l \equiv [-a, -a+d_2 m^{(-)}]$, and for the right region, $\bar{\Omega}^r \equiv [-a+d_2 m^{(+)},a]$.

 We restore the source function $f^\delta(x,y)$ according to the least-squares problem:
\begin{equation}
\label{fdelta}
f^\delta = \mathop{\arg\min}_{\begin{subarray}{c} f\in  C^1(\mathbb{R} \times \bar{\Omega}) \end{subarray} } \sum^{n}_{i=0} \sum_{j\in \mathcal{J}} \left( f(x_i,y_j) -  u^{\delta}_{i,j} (k({u_x})^{\delta}_{i,j} +({u_y})^{\delta}_{i,j}) \right)^2 + \varepsilon \|f\|^2_{H^1(\mathbb{R} \times\Omega)},
\end{equation}
where the regularization parameter satisfies the relation $\varepsilon(\delta)\to 0$ when $\delta\to0$. Numerical experiments indicate that a small value of $\varepsilon$ always produces a satisfactory result. In this paper, we fix $\varepsilon=\delta^2$ in all simulations.

In case we know only the measurements $u^\delta_{i,j}$, we replace the values $\{u^{\delta}_{i,j}, $ $ ({u_x})^{\delta}_{i,j}, $ $ ({u_y})^{\delta}_{i,j} \}$ in \eqref{fdelta} with the smoothed quantities $\{u^\varepsilon, \frac{\partial u^\varepsilon}{\partial x}, \frac{\partial u^\varepsilon}{\partial y}\}$. The function $u^\varepsilon (x,y,t_0)$ is constructed according to the following optimization problem for the left and right regions, respectively:
\begin{multline}
\label{uAlphaL}
u^\varepsilon = \arg\min_{\begin{subarray}{c} v\in C^1(\mathbb{R} \times \bar{\Omega}^l ): \end{subarray} } \frac{1}{n+1}\frac{1}{m^{(-)}+1} \sum^{n}_{i=0}\sum^{m^{(-)}}_{j=0}  \left( v(x_i,y_j,t_0)-u^\delta_{i,j}(t_0) \right)^2 \\
+ \varepsilon^{(-)} (t_0) \left(\left\| \frac{\partial^2 v(x,y,t_0)}{\partial x^2} \right\|^2_{\begin{subarray}{c} L^2(\mathbb{R} \times \bar{\Omega}^l)\end{subarray} }+ \left\| \frac{\partial^2 v(x,y,t_0)}{\partial y^2} \right\|^2_{ \begin{subarray}{c} L^2(\mathbb{R} \times \bar{\Omega}^l) \end{subarray}} \right),
\end{multline}

\begin{multline}
\label{uAlphaR}
u^\varepsilon= \arg\min_{\begin{subarray}{c} v\in C^1(\mathbb{R} \times \bar{\Omega}^r):  \end{subarray} } \frac{1}{n+1}\frac{1}{m-m^{(+)}+1} \sum^{n}_{i=0}\sum^{m}_{j=m^{(+)}}  \left( v(x_i,y_j,t_0)-u^\delta_{i,j}(t_0) \right)^2 \\
+ \varepsilon^{(+)} (t_0) \left(\left\| \frac{\partial^2 v(x,y,t_0)}{\partial x^2} \right\|^2_{\begin{subarray}{c} L^2(\mathbb{R} \times \bar{\Omega}^r) \end{subarray}}+ \left\| \frac{\partial^2 v(x,y,t_0)}{\partial y^2} \right\|^2_{\begin{subarray}{c} L^2(\mathbb{R} \times \bar{\Omega}^r) \end{subarray}} \right),
\end{multline}
where the regularization parameters $\varepsilon^{(\mp)}(t_0)$ satisfy
$$ \displaystyle \frac{1}{n+1}\frac{1}{m^{(-)}+1}  \sum^{n}_{i=0}\sum^{m^{(-)}}_{j=0} \left( u^\varepsilon(x_i,y_j,t_0)-u^\delta_{i,j}(t_0) \right)^2 = \delta^4,$$
$$ \displaystyle \frac{1}{n+1}\frac{1}{m-m^{(+)}+1}  \sum^{n}_{i=0}\sum^{m}_{j=m^{(+)}} \left( u^\varepsilon(x_i,y_j,t_0)-u^\delta_{i,j}(t_0) \right)^2 = \delta^4.$$

According to \cite[Theorem 3.3]{WANG2005121}, the following assertion holds.
\begin{proposition}
\label{NoisyErr}
Suppose that, for a.e. $t\in\bar{\mathcal{T}}$, $u(\cdot,\cdot,t) \in H^2(\mathbb{R} \times \Omega)$. Let $u^\varepsilon(x,y,t)$ be the minimizer of the problems \eqref{uAlphaL} and \eqref{uAlphaR}, with $t_0$ replaced with $t$. Then, for a.e. $t\in\bar{\mathcal{T}}$ and $\varepsilon = \delta^2$,
\begin{align}
\label{NoisyErrIneq}
  \|u^\varepsilon(\cdot,\cdot,t) - u(\cdot,\cdot,t)\|_{H^1(\mathbb{R} \times \Omega)}  \leq C_2 d^{1/4}  + C_3 \sqrt{\delta}, \nonumber
\end{align}
where $d =\max \{ d_1,d_2 \}$, $C_2$ and $C_3$ are constants which depend on the domain $\mathbb{R} \times \Omega$ and the quantity $\left\| \Delta u(x,y,t)\right\|_{L^2(\mathbb{R} \times \Omega)}$.

\end{proposition}

Using Propositions \ref{ProAsympErr} and \ref{NoisyErr} and the triangle inequality $\|f^\delta - $ $f^*\|_{L^2(\mathbb{R} \times \Omega)} \leq \|f^\delta - f_0\|_{L^2(\mathbb{R} \times \Omega)} + \|f_0 - f^*\|_{L^2(\mathbb{R} \times \Omega)}$, it is relatively straightforward to prove the following theorem.
\begin{theorem}
\label{ErrSource}
$f^\delta$, defined in \eqref{fdelta}, is a stable approximation of the exact source function $f^*$ for problem (\textbf{IP}). Moreover, it has the convergence rate
\begin{equation}
\label{ErrSourceIneq}
\|f^\delta - f^*\|_{L^2(\mathbb{R} \times \Omega)} = \mathcal{O} ( \mu \lvert \ln \mu \rvert + d^{1/4} + \sqrt{\delta} ).
\end{equation}
In addition, if $\mu=\mathcal{O}(\delta^{\epsilon + 1/2})$ ($\epsilon$ is any positive number) and $d^{1/4}=\mathcal{O}(\sqrt{\delta})$, the following estimate holds:
\begin{equation*}
\|f^* - f^\delta\|_{L^2(\mathbb{R} \times \Omega)} = \mathcal{O}(\sqrt{\delta}) .
\end{equation*}
\end{theorem}

Based on the above analysis, we build an efficient regularization algorithm for the two-dimensional nonlinear-source inverse problem (\textbf{IP}), as shown below.

\begin{algorithm}[htb]
\caption{Asymptotic expansion-regularization (AER) algorithm for (\textbf{IP}).}
\label{alg:Framwork}
\begin{algorithmic}[1]
\If{The data for the full measurements $\{u^\delta_{i,j}, ({u_x})^\delta_{i,j}, ({u_y})^\delta_{i,j}\}$ are given}
    \State{break;}
\Else
    \If{ Only the measurements $\{ u^\delta_{i,j} \}$ are provided}
        \State Construct the smoothed data $\{u^\varepsilon,\frac{\partial u^\varepsilon}{\partial x},\frac{\partial u^\varepsilon}{\partial y} \}$ by solving \eqref{uAlphaL} and \eqref{uAlphaR};
    \EndIf
\EndIf
\State{Calculate the approximate source function $f^\delta$  using formula \eqref{fdelta}.}
\end{algorithmic}
\end{algorithm}

\section{Proofs of main results} \label{derivationAndProofs}

\subsection{Proof of Theorem \ref{MainThm}}

To prove Theorem \ref{MainThm} and estimate its accuracy (\eqref{eq002}–\eqref{0orderEstim2}), we use the asymptotic method of inequalities \cite{b7}. According to this method, a solution to \eqref{mainproblem} exists if there exist continuous functions $\alpha(x, y, t, \mu)$ and $\beta(x, y, t, \mu)$, called, respectively, lower and upper solutions of \eqref{mainproblem}. First, we recall the definitions of upper and lower solutions and their role in the construction of solution  \eqref{mainproblem} \cite{b7,b8,NEFEDOV201390}.

\begin{definition} \label{Lemma1}
 The functions  $   \beta (x,y, t,\mu) $ and  $ \alpha (x,y,t, \mu)$ are called upper and lower solutions of the problem \eqref{mainproblem}, if they are continuous, twice continuously differentiable in $x$ and $y$, continuously differentiable in $t$, and, for sufficiently small $\mu$, satisfy the following conditions:

  \begin{itemize}[leftmargin=1cm]
 \item[\textbf{(C1):}] $\alpha(x,y,t,\mu)\leq \beta(x,y,t,\mu)$ for $(x,y,t)\in \mathbb{R} \times \bar{\Omega}\times \bar{\mathcal{T}}.$
  \item[\textbf{(C2):}]  $  \displaystyle L[ \alpha]:=\mu \left( \frac{\partial^{2}\alpha}{\partial x^{2}}+\frac{\partial^{2}\alpha}{\partial y^{2}} \right)-\frac{\partial\alpha}{\partial t}+\alpha \left( k\frac{\partial\alpha}{\partial x}+\frac{\partial\alpha}{\partial y}\right)-f(x,y)\geq 0, \quad (x,y,t)\in \mathbb{R} \times \bar{\Omega}\times \bar{\mathcal{T}};$
  \item[\qquad \ \ ] $ \displaystyle L[\beta]:=\mu\left( \frac{\partial^{2}\beta}{\partial x^{2}}+\frac{\partial^{2}\beta}{\partial y^{2}} \right)-\frac{\partial\beta}{\partial t}+\beta \left( k\frac{\partial\beta}{\partial x}+\frac{\partial\beta}{\partial y}\right)-f(x,y )\leq 0, \quad (x,y,t)\in \mathbb{R} \times \bar{\Omega}\times \bar{\mathcal{T}}.$
 \item[\textbf{(C3):}] $\alpha (x,-a,t,\mu)\leq u^{-a} (x)\leq \beta(x,-a,t,\mu), \quad \alpha(x,a,t,\mu)\leq u^{a} (x)\leq \beta(x,a,t,\mu)$.
 \item[\qquad \ \ ] $\alpha (x,y,t,\mu) =\alpha (x+L,y,t,\mu), \quad \beta (x,y,t,\mu) =\beta (x+L,y,t,\mu).$
\end{itemize}
\end{definition}

\begin{lemma} \label{Lemma2}
(\cite{b8}) Let there be an upper $   \beta (x,y,t,\mu) $ and a lower $ \alpha (x,y,t, \mu)$ solution to problem \eqref{mainproblem} satisfying conditions (C1)–(C3) in Definition \ref{Lemma1}. Then, under Assumptions \ref{A1}–\ref{A4}, there exists a solution $u(x,y,t,\mu)$ to problem \eqref{mainproblem} that satisfies the inequalities
$$
\alpha(x,y,t,\mu)\leq u(x,y,t,\mu)\leq \beta(x,y,t,\mu),\ (x,y,t)\in \mathbb{R} \times \bar{\Omega}\times \bar{\mathcal{T}}.
$$
Moreover, the functions $\beta(x,y, t, \mu)$ and $\alpha (x,y, t, \mu)$ satisfy the following estimates:
\begin{align}
& \beta(x,y,t,\mu)-\alpha(x,y,t,\mu) =\mathcal{O}(\mu^{n}), \label{BetaMinusAlpha} \\
& u(x,y,t,\mu)=\alpha(x,y,t,\mu)+\mathcal{O}(\mu^{n})=U_{n-1}(x,y,t,\mu)+\mathcal{O}(\mu^{n}). \label{UMinusAlpha}
\end{align}
 \end{lemma}

\begin{lemma} \label{Lemma3}
(\cite{b7,NEFEDOV201390}) Lemma \ref{Lemma2} also remains valid in the case in which the functions $\alpha(x,y,t,\mu)$ and $\beta(x,y,t,\mu)$ are continuous and their derivatives with respect to $x,y$ have discontinuities  from the class $C^2$  in the direction of the normal to curves on which these solutions are not smooth, and the limit values of the derivatives on the curve $h(x,t)$ satisfy the following condition:
\begin{itemize}[leftmargin=1cm]
\item[\textbf{(C4):}] $\displaystyle \frac{\partial\alpha}{\partial n}(x,h_{\alpha}(x,t)+0,t,\mu)-\frac{\partial\alpha}{\partial n}(x,h_{\alpha}(x,t)-0,t,\mu)\geq 0,$
\end{itemize}
where $h_{\alpha}(x,t)$ is the curve on which the lower solution is not smooth;
\begin{itemize}
 \item[\qquad \   ]  $\displaystyle \frac{\partial\beta}{\partial n}(x,h_{\beta}(x,t)-0,t,\mu)-\frac{\partial\beta}{\partial n}(x,h_{\beta}(x,t)+0,t,\mu)\geq 0,$
 \end{itemize}
 where $h_{\beta}(x,t)$ is the curve on which the upper solution is not smooth.\\
\end{lemma}

The proofs of Lemmas \ref{Lemma2}–\ref{Lemma3} can be found in \cite{b7,b8}. Thus, to prove Theorem \ref{MainThm}, it is necessary to construct the lower and upper solutions  $ \alpha (x,y,t, \mu)$  and $   \beta (x,y,t,\mu) $. Under conditions (C1)–(C4) for  $ \alpha (x,y,t, \mu)$  and $   \beta (x,y,t,\mu) $, estimates \eqref{NorderEstim1} and \eqref{NorderEstim2} will follow directly from Lemma \ref{Lemma2}. Estimate \eqref{NorderEstim3} can be obtained by solving equation \eqref{proofzn2.20} for the quantity $z_{n}(x, y,t, \mu):= u(x, y,t, \mu)-U_{n}(x,y, t, \mu)$, using a Green's function.

We now begin the proof of Theorem \ref{MainThm}.

\begin{proof}
Following the main idea in \cite{b6},
we construct the upper and lower solutions $\beta^{(-)}$, $\beta^{(+)}$, $\alpha^{(-)}$, $\alpha^{(+)}$  and curves $h_{ \beta}$,  $h_{ \alpha}$ as a modification of asymptotic representation \eqref{asymptoticnorder}.

We introduce a positive function $ \rho (x,t) $, which will be defined later in \eqref{phoequat}, and use the notations  $ \rho_{\beta} (x,t)=-\rho (x,t) $ and $ \rho_{\alpha} (x,t)=\rho (x,t) $ to aid in defining the curves $h_{ \beta}(x, t)$ and $h_{ \alpha}(x, t)$, which in turn will determine the position of the inner transition layer for the upper and lower solution, in the form
\begin{align}\label{curveh}
\displaystyle h_{ \beta}(x, t)= \sum_{i=0}^{n+1} \mu^i h_{i}(x,t)+\mu^{n+1} \rho_{\beta}(x, t),  \quad \displaystyle h_{ \alpha}(x, t)= \sum_{i=0}^{n+1} \mu^i h_{i}(x,t)+\mu^{n+1} \rho_{\alpha}(x, t).
\end{align}

In the vicinity of the curve $h_{\beta}(x, t)$,  we pass to the local coordinates $(l, r_{\beta})$  according to the following equations:
$$
x=l-r_{\beta}\sin\alpha_{\beta},\ y=h_{\beta}(l, t)+r_{\beta}\cos\alpha_{\beta}=\hat{h}_{n+1}(l, t)+r_{\beta}\cos\alpha_{\beta}+\mu^{n+1}\rho_{\beta}(l, t),
$$
where $r_{\beta}$ is the distance from the curve $h_{\beta}(x, t)$ along the normal to it, $l$ is the coordinate of the point on the $x$ axis from which this normal is drawn, $\cos{\alpha_{\beta}}=\frac{1}{\sqrt{1+(h_{\beta})_{x}^{2}}}, \ \sin{\alpha_{\beta}}=\frac{(h_{ \beta})_{x}}{\sqrt{1+(h_{\beta})_{x}^{2}}}$, and the derivatives of the function $h_{\beta}$ at each time $t$ are taken at $x=l.$

Similarly, in the vicinity of the curve $h_{\alpha}(x, t)$  we pass to the local coordinates $(l, r_{\alpha})$:
$$
x=l-r_{\alpha}\sin\alpha_{\alpha},\ y=h_{\alpha}(l, t)+r_{\alpha}\cos\alpha_{\alpha}=\hat{h}_{n+1}(l, t)+r_{\alpha}\cos\alpha_{\alpha}+\mu^{n+1}\rho_{\alpha}(l, t),
$$
where $r_{\alpha}$ is the distance from the curve $h_{\alpha}(x, t)$ along the normal to it,  $\cos{\alpha_{\alpha}}=\frac{1}{\sqrt{1+(h_{\alpha})_{x}^{2}}}$, and $ \sin{\alpha_{\alpha}}=\frac{(h_{ \alpha})_{x}}{\sqrt{1+(h_{\alpha})_{x}^{2}}}$.

In the neighborhood of the curves  $h_{\beta}(x, t)$ and  $h_{\alpha}(x, t)$,  we introduce the extended variables
\begin{equation} \label{UpperLowerXi}
\xi_{\beta}=\frac{r_{\beta}}{\mu}, \quad \xi_{\alpha}=\frac{r_{\alpha}}{\mu}.
\end{equation}

The upper and lower solutions of problem \eqref{mainproblem} will be constructed separately in the domains $ \bar{D}_{ \beta}^{(-)}, \bar{D}_{ \beta}^{(+)}$ and $ \bar{D}_{ \alpha}^{(-)}, \bar{D}_{ \alpha}^{(+)}$, in which the surfaces $h_{ \beta}(x, t)$ and $h_{ \alpha}(x, t)$ divide the domain $ \mathbb{R}\times \bar{\Omega}\times \bar{\mathcal{T}} $:
\begin{align}
\label{beta}
\beta(x, y,t, \mu)= \begin{cases}
\beta^{(-)}(x,y,t,\mu ), \ (x,y,t)\in \bar{D}_{ \beta}^{(-)}:=\{(x,y,t) \in \mathbb{R}\times [-a, h_{\beta}] \times  \bar{\mathcal{T}} \},  \\
\beta^{(+)}(x,y,t,\mu ), \ (x,y,t)\in \bar{D}_{ \beta}^{(+)}:=\{(x,y,t) \in \mathbb{R} \times [h_{\beta}, a] \times \bar{\mathcal{T}} \},
\end{cases}
\end{align}
\begin{align}
\label{alpha}
\alpha(x,y, t, \mu)= \begin{cases}
\alpha^{(-)}(x,y,t,\mu ), \ (x,y,t)\in \bar{D}_{ \alpha}^{(-)}:=\{(x,y,t)  \in \mathbb{R}\times [-a, h_{ \alpha}] \times \bar{\mathcal{T}} \}, \\
\alpha^{(+)}(x,y,t,\mu ), \ (x,y,t)\in\bar{D}_{ \alpha}^{(+)}:=\{(x,y,t)  \in \mathbb{R}\times [h_{ \alpha}, a]  \times \bar{\mathcal{T}} \}.
\end{cases}
\end{align}

We match the functions $\beta^{(-)}(x,y, t, \mu)$, $ \beta^{(+)}(x,y, t, \mu)$ and $\alpha^{(-)}(x, y,t, \mu)$, $ \alpha^{(+)}(x,$ $ y,t, \mu)$ on the curves $h_{ \beta}(x, t)$ and $h_{ \alpha}(x, t)$, respectively, so that $\beta(x, t, \mu)$ and $\alpha(x, t, \mu)$ are continuous on these curves and the following equations hold:
\begin{align}
\begin{split} \label{sewingeq}
\displaystyle \beta^{( - )}(x, h_{\beta}, t, \mu)=\beta^{( + )}(x, h_{\beta}, t, \mu)=\frac{\varphi^{(-)}(x, h_{\beta})+\varphi^{(+)}(x, h_{\beta})}{2}, \\
\displaystyle \alpha^{( - )}(x, h_{\alpha}, t, \mu)=\alpha^{( + )}(x, h_{\alpha}, t, \mu)=\frac{\varphi^{(-)}(x, h_{\alpha})+\varphi^{(+)}(x, h_{\alpha})}{2}.
\end{split}
\end{align}

Note that we do not match the derivatives of the upper and lower solutions on the curves $h_{ \beta}(x, t)$ and $h_{ \alpha}(x, t)$, respectively, and so the derivatives $\partial \beta / \partial y$  and $\partial \alpha / \partial y $ have discontinuity points, and therefore fulfillment of condition (C4) is required to hold.

We construct the functions $\beta^{( \mp)}$ and $\alpha^{( \mp)}$ in the following forms:
\begin{align}\label{beta2}
\begin{split}
\beta^{( \mp)}= U_{n+1}^{( \mp)}|_{\xi_{\beta},h_{\beta}}+\mu^{n+1} \left(\epsilon^{( \mp)}(x,y)+q_{0}^{( \mp)}(\xi_{\beta}, t)+\mu q_{1}^{( \mp)}(\xi_{\beta},t) \right), \\
\alpha^{( \mp)}=U_{n+1}^{( \mp)}|_{\xi_{\alpha},h_{\alpha}}-\mu^{n+1} \left( \epsilon^{( \mp)}(x,y)+q_{0}^{( \mp)}(\xi_{\alpha},t)+\mu q_{1}^{( \mp)}(\xi_{\alpha}, t)\right),
\end{split}
\end{align}
where the functions $\epsilon^{( \mp)}(x,y)$ should be designed in such a way that condition (C2) is satisfied for $\beta^{( \mp)}$ and $\alpha^{( \mp)}$ in \eqref{beta2}. The functions $q_{0}^{( \mp)}(\xi_{\beta},t) $ eliminate residuals of order $\mu^n$ arising  in $L[\beta]$ and $L[\alpha]$ and the residuals of order $\mu^{n+1}$  under the condition of continuous matching of the upper solution \eqref{sewingeq}, which arise  as a result of modifying the outer part by adding  $\epsilon^{( \mp)}(x,y)$. The functions $q_{1}^{( \mp)}(\xi_{\beta}, t) $ eliminate  residuals of order $\mu^{n+1}$  arising in $L[\beta]$ by adding   $\epsilon^{( \mp)}(x,y)$ and $q_{0}^{( \mp)} (\xi_{\beta}, t) $.

We now define the functions $\epsilon^{( \mp)}(x,y)$  from the following equations:
\begin{align} \label{epsiloneq}
\begin{split}
\left( k \frac{\partial \epsilon^{( \mp)} }{\partial x} +\frac{\partial \epsilon^{( \mp)} }{\partial y}\right) \varphi^{( \mp)}  +\epsilon^{( \mp)}  \left( k \frac{\partial \varphi^{( \mp)} }{\partial x}+\frac{\partial \varphi^{( \mp)}}{\partial y} \right) = -R  ,\\
\epsilon^{(-)}(x,-a)=R^{(-)}, \quad \epsilon^{(-)}(x,y)=\epsilon^{(-)}(x+L,y),\\ 
\epsilon^{(+)} (x,a)= R^{(+)}, \quad \epsilon^{(+)}(x,y)=\epsilon^{(+)}(x+L,y),
\end{split}
\end{align}
where $R, R^{(-)}, R^{(+)}$ are some positive values independent of $x,y,z$. The functions $\epsilon^{( \mp)} (x,y)$ can be determined explicitly:
\begin{align} \label{epsiloneqexplicit}
\begin{split}
&\epsilon^{(-)} (x,y)= \frac{-R(a+y)+R^{(-)} \varphi^{(-)}(-ak+x-ky,-a)}{\varphi^{(-)} (x,y)},\\
&\epsilon^{(+)} (x,y)= \frac{R(a-y)+R^{(+)} \varphi^{(+)}(ak+x-ky,a)}{\varphi^{(+)} (x,y)},
\end{split}
\end{align}
since $\varphi^{(-)} (x,y) <0 $ and $\varphi^{(+)} (x,y) >0 $, $\epsilon^{( \mp)} (x,y) > 0 $ for $(x,y) \in \mathbb{R} \times \bar{\Omega}$.

 We define the functions $q_{0}^{( \mp)}(\xi_{\beta},t) $ as solutions of the equations
\begin{multline} \label{q0}
\frac{\partial^2 q_{0}^{( \mp)}}{\partial \xi_{\beta}^2}+   \frac{\partial}{\partial \xi_{\beta}} \left(q_{0}^{(\mp)} \frac{{h_0}_{t} ({h_0}_{x}^{2} +1  ) + \tilde{u}(1-k{h_0}_{x} ) }{ \sqrt{1+{h_0}_{x}^{2} }}   \right)=b_{1}^{(\mp)}(\xi_{\beta},l,t) \frac{\partial \rho_{\beta} }{\partial l} \\ 
+b_{2}^{(\mp)}(\xi_{\beta},l,t) \rho_{\beta}  +b_{3}^{(\mp)}(\xi_{\beta},l,t)\frac{\partial \rho_{\beta} }{\partial t}+b_{4}^{(\mp)}(\xi_{\beta},l,t):= H_{q0}^{( \mp)}(\xi_{\beta},  t), 
\end{multline}
where the derivative of the function $ \rho_{\beta} (x, t) $ is taken for $ x = l $ and $b_{1}^{(\mp)}(\xi_{\beta},l,t),$ $b_{2}^{(\mp)}(\xi_{\beta},l,t),$ $b_{3}^{(\mp)}(\xi_{\beta},l,t),$ $b_{4}^{(\mp)}(\xi_{\beta},l,t)$ are known functions, particularly, $b_{3}^{(\mp)}(\xi_{\beta},l,t)=$ $ \displaystyle \frac{\partial Q_{0}^{(\mp)} (\xi_{\beta},l, h_0,t)}{\partial \xi_{\beta}} $ $\sqrt{1+{h_0}_{l}^{2}}.$ 

The boundary conditions for $q_{0}^{( \mp)}(\xi_{\beta},t)$ follow from equation \eqref{sewingeq}, matching the upper solution and taking into account conditions at  $\xi_{\beta} = 0$ for the functions $Q_{i}^{( \mp)}(\xi_{\beta},l,h_\beta, t) $:
\begin{align}  \label{q0conditions}
\begin{split}
 q_{0}^{( \mp)}(0,t)=\frac{\epsilon^{( \pm)}(x,h_0)-\epsilon^{( \mp)}(x,h_0)}{2} \equiv p_{q0}^{( \mp)}(x,h_0), \quad q_{0}^{(\mp)}(\mp\infty,t) = 0.
 \end{split}
\end{align}

Using \eqref{q0} and \eqref{q0conditions}, we write the functions $q_{0}^{(-)}(\xi_{\beta}, t)$ in the explicit form:
\begin{equation} \label{q0explisit}
q_{0}^{( \mp)} =J^{( \mp)}(\xi_{\beta},h_0) \left( p_{q0}^{( \mp)}(x,h_0)   +\int_{0}^{\xi_{\beta}} \frac{1}{J^{( \mp)}(s,h_0)} \int_{\mp \infty}^{s} H_{q0}^{( \mp)} (\eta,t) d\eta ds \right).
\end{equation}

 We define the functions $q_{1}^{( \mp)}(\xi_{\beta}, t) $ from the following equation:
\begin{align} \label{q1}
\begin{split}
\frac{\partial^2 q_{1}^{( \mp)}}{\partial \xi_{\beta}^2}+   \frac{\partial}{\partial \xi_{\beta}} \left(q_{1}^{(\mp)} \frac{{h_0}_{t}  ({h_0}_{x}^{2} +1  ) + \tilde{u}(1-k{h_0}_{x} ) }{ \sqrt{1+{h_0}_{x}^{2} }}   \right)= H_{q1}^{( \mp)}(\xi_{\beta},  t),
\end{split}
\end{align}
where $H_{q1}^{( \mp)}(\xi_{\beta}, t) $ depend on known functions $ h_{0,1}$, $u_{0,1}^{\mp}(x,h_0) $, $Q_{0,1}^{\mp}(\xi_{\beta},l,h_\beta,t)$,  $\rho_{\beta}$, $\epsilon^{( \mp)}(x,h_0)$, and $ q_{0}^{( \mp)}(\xi_{\beta}, t)$. For equation \eqref{q1}, we infer that the boundary condition is
\begin{align*}
q_{1}^{( \mp)}(0,t)=0, \ q_{1}^{( \mp)}(\xi_{\beta},t) \rightarrow 0 \ \text{for} \  \xi_{\beta} \rightarrow \mp \infty.
\end{align*}

Replacing $\rho_{\beta}$ with $\rho_{\alpha}$ and $\xi_{\beta}$ with $\xi_{\alpha}$ in \eqref{q0}–\eqref{q1}, we define the functions $q_{0}^{( \mp)}(\xi_{\alpha}, t)$ and $q_{1}^{( \mp)}(\xi_{\alpha}, t)$ that appear in the functions $\alpha^{( \mp)}$.

The functions $q_{0}^{( \mp)}$ and $q_{1}^{( \mp)}$ satisfy exponential estimates of the types in \eqref{equat22} and \eqref{equat23}.

Now, we have to show that the functions $\beta(x,t,\mu)$ and $\alpha(x,t,\mu)$ are upper and lower solutions to problem \eqref{mainproblem}. To do this, we check conditions (C1)–(C4).

First, we check that condition (C1) has been fulfilled, with regard to the ordering of the lower and upper solutions. To do this, we consider three regions that illustrate the difference between the upper and lower solutions, $\beta-\alpha $:
\begin{align}
\quad \beta-\alpha=\left\{\begin{array}{l}
\beta^{(-)}-\alpha^{(-)}, \quad \RomanNumeralCaps{1} =  \{ (x,y,t):   \mathbb{R} \times [-a, h_{\beta}(x, t)] \times \bar{\mathcal{T}} \},\\
\beta^{(+)}-\alpha^{(-)} , \quad \RomanNumeralCaps{2} = \{ (x,y,t):  \mathbb{R} \times [ h_{\beta}(x, t), h_{\alpha}(x, t)] \times \bar{\mathcal{T}} \},\\
\beta^{(+)}-\alpha^{(+)} , \quad \RomanNumeralCaps{3} = \{ (x,y,t):  \mathbb{R} \times [h_{\alpha}(x, t), a] \times \bar{\mathcal{T}} \}.
\end{array}\right.
\end{align}
First, we find the variables on which $\beta$ and $\alpha $ depend:
\begin{multline*}
y=h_{\beta}(l, t)+r_{\beta}\cos\alpha_{\beta}=h_{\alpha}(l, t)+r_{\alpha}\cos\alpha_{\alpha}= \\
=\hat{h}_{n+1}(l, t)-\mu^{n+1}\rho(l, t)+r_{\beta}\cos\alpha_{\beta}=\hat{h}_{n+1}(l, t)+\mu^{n+1}\rho(l, t)+r_{\alpha}\cos\alpha_{\alpha},
\end{multline*}
$$
\cos\alpha_{\alpha}=\frac{1}{\sqrt{1+(\hat{h}_{n+1})_{x}^{2}}}+\mathcal{O}(\mu^{n+1}), \quad \cos\alpha_{\beta}=\frac{1}{\sqrt{1+(\hat{h}_{n+1})_{x}^{2}}}+\mathcal{O}(\mu^{n+1}),
$$
which we obtain using \eqref{curveh} and \eqref{sincosalpha}.
We rewrite \eqref{xiequation} in the following form:
$$
\xi_{\beta}=\frac{y-h_{\beta}}{\mu}\sqrt{1+({h_{\beta}})_{x}^2} , \quad \xi_{\alpha}=\frac{y-h_{\alpha}}{\mu}\sqrt{1+({h_{\alpha}})_{x}^2},
$$
from which we find
$$
\Delta \xi=\xi_{\beta}-\xi_{\alpha}=2\mu^{n} \left( \rho\sqrt{1+(h_{0})_{x}^{2}}+\frac{\rho_x ({h_0})_{x} (h_0-y)}{\sqrt{1+(h_{0})_{x}^{2}}}\right)+\mathcal{O}(\mu^{n+1}).
$$

For region \RomanNumeralCaps{2}, the following holds:
\begin{align*}
    \begin{split}
       & 0\leq\xi_{\beta}\leq 2\mu^{n}\rho(l, t)\sqrt{1+(h_{0})_{x}^{2}}+\mathcal{O}(\mu^{n+1}),\\
       & 0 \geq\xi_{\alpha}\geq  -2\mu^{n}\rho(l, t)\sqrt{1+(h_{0})_{x}^{2}}+\mathcal{O}(\mu^{n+1}),\\
       & \xi_{\beta}-\xi_{\alpha}=2\mu^{n}  \rho(l, t)\sqrt{1+(h_{0})_{x}^{2}}+\mathcal{O}(\mu^{n+1}).
    \end{split}
\end{align*}

We can write an expression for the difference between the upper and lower solutions:
\begin{align} \label{diffbetaalpha1}
\begin{split}
\beta^{(+)}-\alpha^{(-)} =\sum_{i=0}^{n}\mu^{i} \left(\bar{u}_{i}^{(+)}\left(l-\mu\xi_{\beta}\sin\alpha_{\beta}, h_{\beta}+\mu\xi_{\beta}\cos\alpha_{\beta}\right)+Q_{i}^{(+)}\left(\xi_{\beta}, l, h_{\beta}, t \right) \right)- \\
-\sum_{i=0}^{n}\mu^{i}\left(\bar{u}_{i}^{(-)}\left(l-\mu\xi_{\alpha}\sin\alpha_{\alpha}, h_{\alpha}+\mu\xi_{\alpha}\cos\alpha_{\alpha}\right)+Q_{i}^{(-)} \left(\xi_{\alpha}, l, h_{\alpha}, t \right) \right)+\mathcal{O}(\mu^{n+1}).
\end{split}
\end{align}

Expanding equation \eqref{diffbetaalpha1} in series, and taking into account the notation \eqref{derivativetildeu} and  equation \eqref{sewindcondexpanded0} and the fact that in region \RomanNumeralCaps{2} $\xi_{\beta}=\mathcal{O}(\mu^{n})$ and $ \xi_{\alpha}=\mathcal{O}(\mu^{n})$, we obtain an expression for the difference between the upper and lower solutions in region \RomanNumeralCaps{2}:
\begin{align}
\begin{split}
\beta^{(+)}-\alpha^{(-)} =\frac{\partial Q_{0}^{(+)}}{\partial \xi}(0, l, h_{0}(l, t),t)\xi_{\beta}-\frac{\partial Q_{0}^{(-)}}{\partial \xi}(0, l, h_{0}(l, t),t)\xi_{\alpha}+\mathcal{O}(\mu^{n+1})\\
=2\mu^{n}\rho(l, t)\sqrt{1+(h_{0})_{x}^{2}}\cdot\Phi^{(+)}(0, h_{0}(l, t))+\mathcal{O}(\mu^{n+1}).
\end{split}
\end{align}

 Using equation \eqref{Q0equation} and Assumption \ref{A3}, we can verify that $\Phi^{(+)}(0, h_{0}(l, t))>0$, and, for positive values of $\rho$ and for a sufficiently small $\mu$, we obtain $\beta-\alpha>0, \ (x,y,t)\in \RomanNumeralCaps{2}.$

We now consider the difference between the upper and lower solutions at region \RomanNumeralCaps{3}, where $\xi_{\alpha}\geq 0, \ \xi_{\beta}=\xi_{\alpha}+\Delta \xi. $ Using exponential properties of the functions $Q_{i}^{(+)}$ and $q_{0}^{(+)}$, we obtain
\begin{multline}
\displaystyle \beta-\alpha=\beta^{(+)}-\alpha^{(+)}=\sum_{i=0}^{n}\mu^{i}\left(Q_{i}^{(+)}(\xi_{\beta}, l, h_{\beta}, t)-Q_{i}^{(+)}(\xi_{\alpha}, l, h_{\alpha}, t)\right)\\
+2\mu^{n+1}\epsilon^{(+)}+\displaystyle \mu^{n+1}\left(q_{0}^{(+)}(\xi_{\beta}, t)-q_{0}^{(+)}(\xi_{\alpha}, t)\right)+\mathcal{O}(\mu^{n+2})=2\mu^{n+1}\epsilon^{(+)}\\
+\frac{\partial Q_{0}^{(+)}}{\partial\xi}(\xi_{\alpha}, l, h_{0}, t)(\xi_{\beta}-\xi_{\alpha})+\mathcal{O}(\mu^{n+1})\exp(-\kappa_{1}\xi_{\alpha})+\mathcal{O}(\mu^{n+2}),
\end{multline}
where $\kappa_{1}>0$ is some constant independent of $\xi_\alpha,l,t,\mu$.

Taking into account estimates \eqref{equat22}  and \eqref{equat23} and the equality $ \xi_{\beta}-\xi_{\alpha}=\mathcal{O}(\mu^{n}) $, we obtain an expression for the difference between the upper and lower solutions in region \RomanNumeralCaps{3}:
\begin{equation} \label{diffbetaalpha2}
\beta-\alpha\leq 2\mu^{n+1}\epsilon^{(+)}+\{C_{0}\mu^{n}\exp(-\kappa_{0}\xi_{\alpha})-C_{1}\mu^{n+1}\exp(-\kappa_{1}\xi_{\alpha})\}+\mathcal{O}(\mu^{n+2}),
\end{equation}
where $C_{0}>0$ and $C_{1}>0$ are some constants independent of $\xi_\alpha,l,t,\mu$.

If $\kappa_{0}\leq\kappa_{1}$, the expression in the brackets in \eqref{diffbetaalpha2} is positive, since $ C_{0}>C_{1}\mu$ for a sufficiently small $\mu$. Hence, $\beta-\alpha>0$.

Let $\kappa_{0}>\kappa_{1}$. Consider the region $x\in \mathbb{R}, y \in [h_{\alpha},h_{\alpha}+N\mu\cos{\alpha_{\alpha}}], t\in \bar{\mathcal{T}} ,$ where the value $N>0$. In this region, the value of $r_{\alpha}$ changes on the interval $[0, N\mu]$ and the inequality $\exp(-\kappa_{0}\xi_{\alpha})\geq\exp(-\kappa_{0}N)$ is satisfied, and so the bracketed expression in \eqref{diffbetaalpha2} is positive if $\mu$ is small enough due to the component $C_{0}\mu^{n}\exp(-\kappa_{0}\xi_{\alpha})$. Hence,  $\beta-\alpha>0.$

We now choose the number $N$ that is large enough to satisfy the inequality $C_{1}\exp(-\kappa_{1} $ $N)<2\epsilon^{(+)}$. When $ y \in [h_{\alpha}+N\mu\cos{\alpha_{\alpha}}, a]  $, due to the choice of the number $N$ we obtain
$$
2\mu^{n+1}\epsilon^{(+)}-C_{1}\mu^{n+1}\exp(-\kappa_{1}\xi_{\alpha})\geq\mu^{n+1}(2\epsilon^{(+)}-C_{1}\exp(-\kappa_{1}N))>0.
$$

Thus, $\beta(x, y, t, \mu)-\alpha(x, y, t, \mu)>0$ everywhere in region \RomanNumeralCaps{3}. The proof of the inequality $\beta(x, y, t, \mu)-\alpha(x, y, t, \mu)>0$ for region \RomanNumeralCaps{1} is developed in the same way as for region \RomanNumeralCaps{3}.

The method of constructing the upper and lower solutions implies the inequalities
\begin{equation*}
L[\beta]=-\mu^{n+1} R + \mathcal{O}(\mu^{n+2})<0, \quad L[\alpha]=\mu^{n+1} R + \mathcal{O}(\mu^{n+2})>0,
\end{equation*}
where $R$ is a constant from \eqref{epsiloneqexplicit}. This verifies condition (C2).

Condition (C3) is satisfied for sufficiently large values $R^{(-)}$ and $R^{(+)}$ in the boundary conditions of equation \eqref{epsiloneq}.

We now check condition (C4) for the upper solution, and expand it in powers of $\mu$; due to the matching of the formal asymptotics \eqref{sewindcondexpanded0} and \eqref{matchingfirstord}, the coefficients at $\mu^{i}$ for $i = 1,\ldots,n$ are equal to zero, and the coefficient at $\mu^{n + 1}$ includes only terms that arise as a result of the modification of the asymptotic expansion:
\begin{align}
\mu \left( \frac{\partial\beta^{( - )}}{\partial n}-\frac{\partial\beta^{( + )}}{\partial n} \right)\Big|_{h=h_{\beta}} = \mu^{n+1} \left( \frac{\partial {q_0}^{( - )}}{\partial \xi_{\beta}}(0,t) - \frac{\partial q_{0}^{( + )}}{\partial \xi_{\beta}}(0,t) \right) +\mathcal{O}(\mu^{n+2}).
\end{align}

Using the explicit solution for $  {q_0}^{( \mp)}(0,t)$ \eqref{q0explisit}, we find
\begin{multline} \label{q0lminusq0r}
 \frac{\partial {q_0}^{( - )}}{\partial \xi_{\beta}}(0,t) - \frac{\partial q_{0}^{( + )}}{\partial \xi_{\beta}}(0,t) =\frac{\partial \rho_{\beta}(x, t) }{\partial t} \left( \int_{- \infty}^{0} b_{3}^{(-)}(s,x,t) ds + \int_{0}^{+ \infty}b_{3}^{(+)}(s,x,t) ds \right) \\
+\frac{\partial \rho_{\beta}(x, t) }{\partial x} H_1(x,t)+\rho_{\beta} (x, t) H_2(x,t)+H_3(x,t),
\end{multline}
where
\begin{flalign*}
H_1(x,t)=  \int_{- \infty}^{0} b_{1}^{(-)}(  s,x,t) ds + \int_{0}^{+ \infty} b_{1}^{(+)}(s,x,t) ds  ,\\
H_2(x,t)=  \int_{- \infty}^{0} b_{2}^{(-)}(s,x,t) ds + \int_{0}^{+ \infty} b_{2}^{(+)}(s,x,t) ds ,
\end{flalign*}
\begin{flalign*}
H_3(x,t)=\left(\epsilon^{( +)} -\epsilon^{( -)} \right) \left(  \frac{- {h_0}_{t} (x,t) ({h_0}_{x}^{2} (x,t)+1  ) + \varphi (x,h_0(x,t))(k{h_0}_{x} (x,t)-1) }{ \sqrt{1+{h_0}_{x}^{2} (x,t)}} \right) \\
+\int_{- \infty}^{0} b_{4}^{(-)}(s,x,t) ds + \int_{0}^{+ \infty} b_{4}^{(+)}(s,x,t) ds.
\end{flalign*}

We choose the function $\rho_{\beta}(x, t) $ as a solution to the problem
\begin{align} \label{phoequat}
\begin{split}
&\frac{\partial \rho_{\beta} }{\partial t} \left( \varphi^{(+)}-\varphi^{(-)} \right)\sqrt{1+{h_0}_{x}^{2}} = -H_1(x,t) \frac{\partial \rho_{\beta} }{\partial x}-H_2(x,t)\rho_{\beta}-H_3(x,t)+\sigma, \\
&\rho_{\beta} (x,0)=\rho^0(x), \quad  \rho_{\beta} (x,t)=\rho_{\beta} (x+L,t), \quad t \in \bar{\mathcal{T}}.
\end{split}
\end{align}
where $\sigma$ is a constant independent of $x,y,t$.

Since the difference $\varphi^{(+)}-\varphi^{(-)} $ is positive and we set the constant $\sigma$ and the value $\rho^0(x)$ to be positive for any $x$,   the solution $\rho (x,t)$ to equation \eqref{phoequat} is also positive for a sufficiently big $\sigma$.

For such $\rho (x,t)$, we obtain:
\begin{equation}
\mu \left( \frac{\partial\beta^{( - )}}{\partial n}-\frac{\partial\beta^{( + )}}{\partial n} \right)\Big|_{h(l,t)=h_{\beta}(l,t)} =\mu^{n+1} \sigma +\mathcal{O}(\mu^{n+2}) >0.
\end{equation}

Similarly, condition (C4) is satisfied for the functions $\alpha^{( \mp)}$, and the constructed upper and lower solutions guarantee the existence of a solution $u (x, t, \mu)$ to the problem \eqref{mainproblem} satisfying the inequalities
\begin{equation}
 \alpha(x,y,t,\mu)\leq u(x,y,t, \mu)\leq \beta(x,y,t,\mu).
\end{equation}
In addition, estimates \eqref{NorderEstim1} and \eqref{NorderEstim2} are valid.

We now show that estimate \eqref{NorderEstim3} also holds. First, we estimate the difference  $z_{n}(x,y, t,$ $ \mu)\equiv u(x,y, t, \mu)-U_{n}(x,y, t, \mu) $; the function $z_{n}(x,y,t, \mu)$ satisfies the equation
\begin{align} \label{proofzn2.20}
\mu \left( \frac{\partial^{2}z_{n}}{\partial x^{2}} +\frac{\partial^{2}z_{n}}{\partial y^{2}}\right)-\frac{\partial z_{n}}{\partial t}- k\left(  U_{n} \frac{\partial U_{n}}{\partial x} -  u \frac{\partial u}{\partial x}  \right) - \left( U_{n} \frac{\partial U_{n}}{\partial y} -  u \frac{\partial u}{\partial y}  \right)  =\mu^{n+1}\psi(x,y, t, \mu),
\end{align}
for $(x,y, t)\in \mathbb{R} \times \bar{\Omega}\times \bar{\mathcal{T}}$, with zero boundary conditions, where $|\psi(x,y, t, \mu)|\leq c_{0}$ and $c_{0} $ is a constant independent of $x,y, t, \mu $. Using the estimates from Lemma \ref{Lemma2}, we obtain
\begin{equation} \label{proof2.21}
z_{n}(x,y, t, \mu)= u(x,y, t, \mu)-U_{n}(x, y,t, \mu)= \mathcal{O} (\mu^{n+1}).
\end{equation}

The terms of equation \eqref{proofzn2.20} can be represented in the following form:
\begin{equation}
 k U_{n} \frac{\partial U_{n}}{\partial x}- k u \frac{\partial u}{\partial x}=\frac{\partial}{\partial x}\int_{u}^{U_{n}} (k s) ds, \quad   U_{n} \frac{\partial U_{n}}{\partial y}-  u \frac{\partial u}{\partial y}=\frac{\partial}{\partial y}\int_{u}^{U_{n}} (s) ds.
\end{equation}

We then rewrite \eqref{proofzn2.20} in the following form:
\begin{align} \label{proof2.22}
\frac{\partial^{2}z_{n}}{\partial x^{2}}+\frac{\partial^{2}z_{n}}{\partial y^{2}}- \frac{1}{\mu} \frac{\partial z_{n}}{\partial t}-K z_{n}  =-K z_{n}+\frac{1}{\mu}\frac{\partial}{\partial x}\int_{u}^{U_{n}} (k s) ds +\frac{1}{\mu}\frac{\partial}{\partial y}\int_{u}^{U_{n}} ( s) ds+ \mu^{n}\psi(x,y, t, \mu).
\end{align}

We define
$$
\theta (x,y, t, \mu):=\mu^{n}\psi(x,y, t, \mu),
$$
and, changing the variable to $\tilde{t}=\mu t$, we can rewrite \eqref{proof2.22} in the following form:
\begin{align} \label{proof2.22v2}
\frac{\partial^{2}z_{n}}{\partial x^{2}}+\frac{\partial^{2}z_{n}}{\partial y^{2}}- \frac{1}{\mu} \frac{\partial z_{n}}{\partial t}-K z_{n}
=\frac{1}{\mu}\left( \frac{\partial}{\partial x}\int_{u}^{U_{n}} (k s) ds  +\frac{\partial}{\partial y}\int_{u}^{U_{n}} ( s) ds \right)-K z_{n}+ \theta(x,y, \frac{\tilde{t}}{\mu}, \mu).
\end{align}

Using a Green's function for the parabolic operator on the left-hand side of \eqref{proof2.22v2}, for any $ (x,y,t)\in \mathbb{R} \times \bar{\Omega}\times \bar{\mathcal{T}}$,  we obtain the representation for $z_{n}$ \cite{pao1992}:
\begin{multline} \label{proof2.23}
z_{n}=\int_{0}^{L} \int_{-a}^{a} G(x, y, \mu t, \zeta,\zeta_1, \mu t_{0})z_{n}(\zeta,\zeta_1, \mu t_{0})d\zeta d\zeta_1\\
-\int_{\mu t_{0}}^{\mu t}d\tau\int_{0}^{L} \int_{-a}^{a} G(x,y, \mu t, \zeta,\zeta_1, \frac{\tau}{\mu} )
\displaystyle \times \left(-K z_{n}(\zeta,\zeta_1, \frac{\tau}{\mu})+\theta (\zeta,\zeta_1, \frac{\tau}{\mu}, \mu)   \vphantom{\int_{\frac{\tau}{\mu}}^{U_{n}(\zeta,\zeta_1,\frac{\tau}{\mu},\mu)}} \right.\\
\left.+\frac{1}{\mu}\frac{\partial}{\partial\zeta}\int_{u(\zeta,\zeta_1,\frac{\tau}{\mu},\mu)}^{U_{n}(\zeta,\zeta_1,\frac{\tau}{\mu},\mu)}(k s)ds  +\frac{1}{\mu}\frac{\partial}{\partial\zeta_1}\int_{u(\zeta,\zeta_1,\frac{\tau}{\mu},\mu)}^{U_{n}(\zeta,\zeta_1,\frac{\tau}{\mu},\mu)}( s)ds \right)d\zeta d\zeta_1.
\end{multline}

Using integration by parts and the boundary conditions for $G$, we can transform the fourth term in \eqref{proof2.23} as follows:
\begin{multline} \label{proof2.24}
 \int_{\mu t_{0}}^{\mu t}d\tau\int_{0}^{L} \int_{-a}^{a} G(x,y,\mu t, \zeta,\zeta_1, \frac{\tau}{\mu})\frac{1}{\mu}\frac{\partial}{\partial\zeta}\int_{u(\zeta,\zeta_1,\frac{\tau}{\mu},\mu)}^{U_{n}(\zeta,\zeta_1,\frac{\tau}{\mu},\mu)}(k s)dsd\zeta d\zeta_1 \\
 =-\int_{\mu t_{0}}^{\mu t}d\tau\int_{0}^{L} \int_{-a}^{a} G_{\zeta}(x,y, \mu t, \zeta,\zeta_1, \frac{\tau}{\mu})\frac{1}{\mu}\int_{u(\zeta,\zeta_1,\frac{\tau}{\mu},\mu)}^{U_{n}(\zeta,\zeta_1,\frac{\tau}{\mu},\mu)}(k s)dsd\zeta d\zeta_1 \\
=-\int_{\mu t_{0}}^{\mu t}d\tau\int_{0}^{L} \int_{-a}^{a} G_{x}(x,y, \mu t, \zeta,\zeta_1, \frac{\tau}{\mu})\frac{1}{\mu}\int_{u(\zeta,\zeta_1,\frac{\tau}{\mu},\mu)}^{U_{n}(\zeta,\zeta_1,\frac{\tau}{\mu},\mu)}(k s)dsd\zeta d\zeta_1 \\
=-\displaystyle \frac{\partial}{\partial x} \left( \int_{\mu t_{0}}^{\mu t}d\tau\int_{0}^{L} \int_{-a}^{a} G(x,y, \mu t, \zeta,\zeta_1, \frac{\tau}{\mu})\frac{1}{\mu}\int_{u(\zeta,\zeta_1,\frac{\tau}{\mu},\mu)}^{U_{n}(\zeta,\zeta_1,\frac{\tau}{\mu},\mu)}(k s)dsd\zeta d\zeta_1 \right).
\end{multline}
In a similar way, we rewrite the last term in \eqref{proof2.23},  and  from \eqref{proof2.23} we obtain the following representation for the derivative $\displaystyle \frac{\partial z_{n}}{\partial x}$:
\begin{multline} \label{proof2.25}
\frac{\partial z_{n}}{\partial x}=\int_{0}^{L} \int_{-a}^{a} G_{x}(x,y, \mu t, \zeta, \zeta_1, \mu t_{0})z_{n}(\zeta, \zeta_1, \mu t_{0})d\zeta d\zeta_1 \\
-\int_{\mu t_{0}}^{\mu t}d\tau\int_{0}^{L} \int_{-a}^{a} G_{x}(x,y, \mu t, \zeta, \frac{\tau}{\mu})\left(-K z_{n}(\zeta,\zeta_1, \frac{\tau}{\mu})+\theta(\zeta, \zeta_1, \frac{\tau}{\mu}, \mu)\right)d\zeta d\zeta_1\\
+\displaystyle \frac{\partial^{2}}{\partial x^{2}} \left(\int_{\mu t_{0}}^{\mu t}d\tau\int_{0}^{L} \int_{-a}^{a} G(x,y, \mu t, \zeta,\zeta_1, \frac{\tau}{\mu})\frac{1}{\mu}\int_{u(\zeta, \zeta_1,\frac{\tau}{\mu},\mu)}^{U_{n}(\zeta, \zeta_1,\frac{\tau}{\mu},\mu)}(k s)dsd\zeta d\zeta_1 \right) \\
+\displaystyle \frac{\partial^{2}}{\partial x \partial y } \left(\int_{\mu t_{0}}^{\mu t}d\tau\int_{0}^{L} \int_{-a}^{a} G(x,y, \mu t, \zeta,\zeta_1, \frac{\tau}{\mu})\frac{1}{\mu}\int_{u(\zeta, \zeta_1,\frac{\tau}{\mu},\mu)}^{U_{n}(\zeta, \zeta_1,\frac{\tau}{\mu},\mu)}(s)dsd\zeta d\zeta_1 \right).
\end{multline}

The validity of representation \eqref{proof2.25} follows from the estimates  \cite[Page 49]{pao1992}:
\begin{equation*}\label{estimat1}
\left|\int_{0}^{L} \int_{-a}^{a} G_{x}(x,y, \mu t, \zeta,\zeta_1 , \mu t_{0})d\zeta d\zeta_1 \right|\leq C_4,
\end{equation*}
\begin{equation*}\label{estimat2}
 \left|\int_{\mu t_{0}}^{\mu t}d\tau\int_{0}^{L} \int_{-a}^{a} G_{x}(x,y, \mu t, \zeta,\zeta_1, \frac{\tau}{\mu} )d\zeta d\zeta_1 \right|\leq C_4,
\end{equation*}
\begin{equation*}\label{estimat3}
\left|\frac{\partial^{2}}{\partial x^{2}}\int_{\mu t_{0}}^{\mu t}d\tau\int_{0}^{L} \int_{-a}^{a} G(x,y, \mu t, \zeta, \zeta_1, \frac{\tau}{\mu})\frac{1}{\mu}d\zeta d\zeta_1 \right|\leq C_4,
\end{equation*}
\begin{equation*}\label{estimat4}
\left|\frac{\partial^{2}}{\partial x \partial y }\int_{\mu t_{0}}^{\mu t}d\tau\int_{0}^{L} \int_{-a}^{a} G(x,y, \mu t, \zeta, \zeta_1, \frac{\tau}{\mu})\frac{1}{\mu}d\zeta d\zeta_1 \right|\leq C_4,
\end{equation*}
where $C_4$ is a constant independent of $(x,y,t)$ and $(\zeta, \zeta_1, \tau)$.
We determine that the first term of representation \eqref{proof2.25} has the estimate $\mathcal{O}(\mu^{n+1})$, and the second term has the estimate $\mathcal{O}(\mu^{n})$. The last two terms in representation \eqref{proof2.25} can be estimated by
$$
\frac{1}{\mu} \left|\int_{u(\zeta, \zeta_1, \frac{\tau}{\mu},\mu)}^{U_{n}(\zeta,\zeta_1,\frac{\tau}{\mu},\mu)}( s)ds\right| = \mathcal{O}(\mu^{n}) .
$$
Using these estimates, from \eqref{proof2.25} we obtain $\displaystyle \frac{\partial z_{n}}{\partial x}=\mathcal{O}(\mu^{n})$ for $(x,y, t)\in \mathbb{R} \times \bar{\Omega}\times \bar{\mathcal{T}}$. Similarly, we obtain the estimation $\displaystyle \frac{\partial z_{n}}{\partial y}=\mathcal{O}(\mu^{n})$, and the estimation $\displaystyle \frac{\partial z_{n}}{\partial n}=\mathcal{O}(\mu^{n})$ follows directly from \eqref{NormalToCurveEq}. This completes the proof of Theorem \ref{MainThm}.
\end{proof}

\subsection{Proof of Proposition \ref{ProAsympErr}}

\begin{proof}
First, we note that the exact source function $f^*$ has the following representation according to equations \eqref{zeroorderregularequation1} and \eqref{zeroorderregularequation2}:

\begin{align} \label{equatforf*}
f^*= \begin{cases}
\displaystyle \bar{u}_{0}^{(-)} \left(k\frac{\partial\bar{u}_{0}^{(-)}}{\partial x}+\frac{\partial\bar{u}_{0}^{(-)}}{\partial y}\right),  \quad (x,y)\in \mathbb{R} \times  (-a,h_0(x,t) - \Delta h/2), \\ \\
\displaystyle \bar{u}_{0}^{(+)}\left( k\frac{\partial\bar{u}_{0}^{(+)}}{\partial x}+\frac{\partial\bar{u}_{0}^{(+)}}{\partial y}\right), \quad (x,y)\in \mathbb{R} \times  (h_0(x,t) + \Delta h/2,a).
  \end{cases}
\end{align}

Let $\Omega'=(-a,h_0(x,t) - \Delta h/2)$. By using Corollary \ref{Corollary1}, we deduce, together with $\lvert \Omega' \rvert \leq \lvert \Omega \rvert=a$, that

\begin{align} \label{eq004IPPf}
\begin{split}
&\left\| \varphi^{(-)}(x,y) - u (x,y,t) \right\|_{L^p(\mathbb{R} \times\Omega')} \leq C \mu,\\
&\left\| \frac{\partial \varphi^{(-)}}{\partial x} - \frac{\partial  u}{\partial x} \right\|_{L^p(\mathbb{R} \times\Omega')} \leq C \mu , \quad \left\| \frac{\partial  \varphi^{(-)}}{\partial y} - \frac{\partial  u}{\partial y} \right\|_{L^p(\mathbb{R} \times\Omega')} \leq C \mu,\\
&\left\| \frac{\partial  \varphi^{(-)}}{\partial x} \right\|_{L^p(\mathbb{R} \times \Omega')} \leq \left\| \frac{\partial  \varphi^{(-)}}{\partial x} - \frac{\partial  u}{\partial x} \right\|_{L^p(\mathbb{R} \times\Omega')} + \left\|\frac{\partial  u}{\partial x} \right\|_{L^p(\mathbb{R} \times \Omega')}\leq C \mu + \left\|\frac{\partial  u}{\partial x} \right\|_{L^p(\mathbb{R} \times \Omega')},\\
&\left\| \frac{\partial  \varphi^{(-)}}{\partial y} \right\|_{L^p(\mathbb{R} \times \Omega')} \leq \left\| \frac{\partial  \varphi^{(-)}}{\partial y} - \frac{\partial  u}{\partial y} \right\|_{L^p(\mathbb{R} \times \Omega')} + \left\|\frac{\partial  u}{\partial y} \right\|_{L^p(\mathbb{R} \times \Omega')} \leq C \mu + \left\|\frac{\partial  u}{\partial y} \right\|_{L^p(\mathbb{R} \times \Omega')}.
\end{split}
\end{align}

From estimates \eqref{eq004IPPf}, we conclude that
\begin{align} \label{proofPro1}
\begin{split}
& \left\|f^* - f_0 \right\|_{L^p(\mathbb{R} \times \Omega')}  = \left\| \varphi^{(-)} \left( k \frac{\partial \varphi^{(-)}}{\partial x} +\frac{\partial \varphi^{(-)}}{\partial y} \right) - u  \left( k \frac{\partial u }{\partial x} + \frac{\partial u }{\partial y} \right) \right\|_{L^p(\mathbb{R} \times \Omega')} \\
& \quad \leq \left\| \varphi^{(-)} \left( k \frac{\partial \varphi^{(-)}}{\partial x} +\frac{\partial \varphi^{(-)}}{\partial y} \right) - u  \left( k \frac{\partial \varphi^{(-)} }{\partial x} + \frac{\partial \varphi^{(-)} }{\partial y} \right) \right\|_{L^p(\mathbb{R} \times \Omega')}   \\
& \qquad \qquad \qquad \qquad \qquad   +\left\| u \left( k \frac{\partial \varphi^{(-)}}{\partial x} +\frac{\partial \varphi^{(-)}}{\partial y} \right) - u  \left( k \frac{\partial u }{\partial x} + \frac{\partial u }{\partial y} \right) \right\|_{L^p(\mathbb{R} \times \Omega')} \\
 & \quad \leq \left\| \varphi^{(-)} - u  \right\|_{L^p(\mathbb{R} \times \Omega')} \left\| k \frac{\partial \varphi^{(-)}}{\partial x} +\frac{\partial \varphi^{(-)}}{\partial y} \right\|_{L^p(\mathbb{R} \times \Omega')} \\
& \qquad \qquad \qquad \qquad \qquad   +   \left\| u  \right\|_{L^p(\mathbb{R} \times \Omega')} \left\|  k \frac{\partial \varphi^{(-)}}{\partial x} +\frac{\partial \varphi^{(-)}}{\partial y} -k \frac{\partial u }{\partial x} - \frac{\partial u }{\partial y} \right\|_{L^p(\mathbb{R} \times \Omega')}\\
 & \quad   \leq C (k+1) \left( C \mu +\left\|u \right\|_{W^{1,p}(\mathbb{R} \times \Omega')} \right) \mu =: c_1  \mu.
\end{split}
\end{align}

Following exactly the same reasoning process, we also derive the following inequality:
\begin{equation} \label{proofPro2}
\left\|f^* - f_0 \right\|_{L^p(\mathbb{R} \times [h_0 + \Delta h/2,a])}\leq c_2 \mu
\end{equation}
with a constant $c_2$. Now, we note that
\begin{align} \label{proofPro3}
\displaystyle \left\|f^* - f_0 \right\|_{L^p(\mathbb{R} \times [h_0 - \Delta h/2, h_0 + \Delta/2])}\leq
\left\|f^* \right\|_{L^p(\mathbb{R} \times [h_0 - \Delta h/2, h_0 + \Delta h/2])} 
+ \left\| f_0 \right\|_{L^p(\mathbb{R} \times [h_0 - \Delta h/2, h_0 + \Delta h/2])} \leq c_3 \mu \lvert \ln \mu \rvert
\end{align}
with $c_3=\left\|f^* \right\|_{C(\mathbb{R} \times \Omega)} + \left\| f_0 \right\|_{C(\mathbb{R} \times \Omega)}$. By combining \eqref{proofPro1}–\eqref{proofPro3}, we deduce that
\begin{multline*}
\left\|f^* - f_0 \right\|^p_{L^p(\mathbb{R} \times [-a,a])}  =  \left\|f^* - f_0 \right\|^p_{L^p(\mathbb{R} \times [-a,h_0 - \frac{\Delta h}{2}])} + \left\|f^* - f_0 \right\|^p_{L^p(\mathbb{R} \times [h_0 + \frac{\Delta h}{2},a] )} \\
 + \left\|f^* - f_0 \right\|^p_{L^p(\mathbb{R} \times [h_0 - \frac{\Delta h}{2},h_0 + \frac{\Delta h}{2}])}   \leq c^p_1 \mu^p + c^p_2 \mu^p + c^p_3 \mu^p \lvert \ln \mu \rvert^p ,
\end{multline*}
which yields the required estimate \eqref{f0Ineq}, with $C_1 = \left( c^p_1  + c^p_2 + c^p_3 \right)^{1/p}$.
\end{proof}

\subsection{Proof of Theorem \ref{ErrSource}}

\begin{proof}
Let $f_0$ be the pre-approximate source function defined in \eqref{f0}. Similarly to the proof of Proposition \ref{ProAsympErr}, for small enough $\delta$ we obtain the following inequalities:
\begin{equation}
\label{f0falphaErr}
\begin{array}{ll}
\displaystyle  \|f^\delta - f_0  \|_{L^2(\mathbb{R} \times \Omega')}  = \left\| u^\delta \left( k \frac{\partial u^\delta }{\partial x} +\frac{\partial u^\delta }{\partial y} \right) - u  \left( k \frac{\partial u }{\partial x} + \frac{\partial u }{\partial y} \right) \right\|_{L^2(\mathbb{R} \times \Omega')}   \\
 \displaystyle  \quad \leq \left\| u^\delta - u  \right\|_{L^2(\mathbb{R} \times \Omega')} \left\| k \frac{\partial u^\delta}{\partial x} +\frac{\partial u^\delta}{\partial y} \right\|_{L^2(\mathbb{R} \times \Omega')}   \\ 
  \displaystyle \qquad \qquad \qquad \qquad  \qquad \qquad  +   \left\| u  \right\|_{L^2(\mathbb{R} \times \Omega')} \left\|  k \frac{\partial u^\delta}{\partial x} +\frac{\partial u^\delta}{\partial y} -k \frac{\partial u }{\partial x} - \frac{\partial u }{\partial y} \right\|_{L^2(\mathbb{R} \times \Omega')}\\
 \displaystyle  \quad \leq \left\| u^\delta - u  \right\|_{L^2(\mathbb{R} \times \Omega')} (k+1) \left( \|u^\delta - u   \|_{H^{1}(\mathbb{R} \times \Omega')} + \|u   \|_{H^{1}(\mathbb{R} \times \Omega')} \right)  +   \left\| u  \right\|_{L^2(\mathbb{R} \times \Omega')}  (k+1) \|u^\delta - u   \|_{H^{1}(\mathbb{R} \times \Omega')} \\
\displaystyle  \quad
\leq 2(k+1)  \left( \|  u \|_{H^{1}(\mathbb{R} \times \Omega')} + 1 \right)  \|u^\delta - u   \|_{H^{1}(\mathbb{R} \times \Omega')} \\  \displaystyle  \quad
\leq 2(k+1) \left( \|  u \|_{H^{1}(\mathbb{R} \times \Omega')} + 1 \right) \left( C_2 d^{1/4} + C_3 \sqrt{\delta}  \right)  \\ \displaystyle  \quad
\leq C_5 (d^{1/4} + \sqrt{\delta}),
\end{array}
\end{equation}
where $\displaystyle  C_5 = 2(k+1) \left( \|  u \|_{H^{1}(\mathbb{R} \times \Omega')} + 1 \right) \left(   C_2  + C_3  \right)$.

By combining \eqref{f0Ineq} and \eqref{f0falphaErr}, we can derive the estimate \eqref{ErrSourceIneq} from the triangle inequality $\|f^\delta - f^*  \|_{L^2(\mathbb{R} \times\Omega)}\leq \|f^\delta - f_0 \|_{L^2(\mathbb{R} \times\Omega)} + \|f_0 -f^* \|_{L^2(\mathbb{R} \times\Omega)}  $.
\end{proof}

\section{Numerical examples}
\label{simulation}
In this section, some numerical experiments are presented to illustrate the efficiency of our new approach. For each example, we first solve the forward problem, where we are looking for a solution to equation \eqref{mainproblem} in the form of an autowave, the existence and uniqueness of which is guaranteed by Theorem \ref{MainThm}. Under Assumptions \ref{A1}–\ref{A4}, the asymptotic solution to problem \eqref{mainproblem} has the following form:
\begin{align} \label{asymptoticsolEXAMPLE1}
U_{0}=\begin{cases}
\displaystyle \varphi^{(-)} (x,y) +\frac{\varphi^{(+)}(x,h_0)-\varphi^{(-)}(x,h_0)}{\exp \left(\frac{\left(h_0-y \right) \left(\varphi^{(+)}(x,h_0)-\varphi^{(-)}(x,h_0)\right) \left(1-k{h_0}_{x} \right) }{2\mu }  \right)+1}, \  y \in [-2,h_0],\\
\displaystyle \varphi^{(+)} (x,y)+\frac{\varphi^{(-)}(x,h_0)-\varphi^{(+)}(x,h_0)}{\exp \left(\frac{\left(h_0-y \right) \left(\varphi^{(-)}(x,h_0)-\varphi^{(+)}(x,h_0)\right) \left(1-k{h_0}_{x} \right) }{2\mu }  \right)+1}, \  y \in [h_0,2].
\end{cases}
\end{align}
where $\varphi^{(-)} $ and $\varphi^{(+)} $ are the solutions to the reduced equations \eqref{zeroorderregularequation1} and \eqref{zeroorderregularequation2}.

For convenience, we write the equation \eqref{h0mainequation} below and verify Assumption \ref{A3} by solving it:
\begin{align*}
\begin{cases}
\displaystyle {h_0}_{t} \left( 1+{h_0}_{x}^{2}\right)=\left( {h_0}_{x}-\frac{1}{2} \right) \left( \varphi^{(-)} (x,h_0)+\varphi^{(+)} (x,h_0)\right),\\
{h_0}(x, t)={h_0}(x+L, t), \quad h_0(x,0)=h_0^{*}=0.
\end{cases}
\end{align*}

To satisfy Assumption \ref{A4}, we take the initial function in \eqref{mainproblem} as follows:
\begin{equation}
\label{InitialEx}
\displaystyle u_{init} (x,y,\mu )=\frac{u^{a}-u^{-a}}{2} \tanh\left(x+\frac{y-h_0^{*}}{\mu}\right)+ \frac{u^{a}+u^{-a}}{2},
\end{equation}
 with an inner transition layer in the vicinity of $y=h_0^{*}$.

For the simulation of inverse problems \textbf{(IP)} for determining the source function $f(x,y)$,  we first solve the forward problem numerically using the finite-volume method, where we introduce the mesh uniformly with respect to spatial variables $x,y$.

In the examples, we consider the case in which only measurements $\{u^\delta_{i,j}\}^{n,m}_{i,j=0}$ are provided, as a more complex scenario. We take the data $u(x_i,y_j,t_0)$ obtained from the numerical result for the forward problem (Fig. \ref{fig:numericsolexample1})  which belongs to the mesh knots $X_n=\lbrace x_i, \ 0 \leq i \leq n: x_i= x_0+d_1 i, \ d_1 = 1/n \rbrace $ and $Y_m=\lbrace y_j, \ 0 \leq j \leq m: y_m= -a+d_2 j, \ d_2 = 1/m \rbrace $.

We skip the points from transition layer $(h_0(x,t_0)-\Delta h /2, h_0(x,t_0)+\Delta h/2)$, and use nodes in only two regions, located on two sides of the transition layer, with indices  $j= 0, \cdots, m^l$, $j=m^r, \cdots, m$, and $i= 0, \cdots, n$. The uniform noise \eqref{noiseUni} is added to the values $u(x_i, y_j, t_0)$ to generate the artificial noisy data  $\{u^{\delta}_{i,j}\}$ on the left and right intervals with respect to the transition layer:
\begin{equation}
\label{noiseUni}
u^{\delta}_{i,j} := [1+\delta(2\cdot \text{rand} -1)] \cdot u(x_i,y_j,t_0),
\end{equation}
where ``$\text{rand}$'' returns a pseudo-random value drawn from a uniform distribution on $[0,1]$.

Following Algorithm \ref{alg:Framwork}, we obtain the smoothed data $\{u^\varepsilon,\frac{\partial u^\varepsilon}{\partial x},\frac{\partial u^\varepsilon}{\partial y} \}$ in accordance with optimization problems \eqref{uAlphaL} and \eqref{uAlphaR}, for the left and right segments, respectively.

Finally, we calculate the regularized approximate source function $f^\delta(x,y)$ using formula \eqref{fdelta}; the retrieved source function is then compared with the exact source function $f^*$.

\subsection{Example 1}
\subsubsection{Forward problem}

Consider the following reaction–diffusion–advection equation in a two-dimensional setting with a periodic boundary condition along the $ x $ axis and with the given source function $f(x,y)=\cos{ ( \pi x /4)} \cos{(\pi y /4)}$:
\begin{align} \label{forwardexample1}
\begin{cases}
\displaystyle 0.08 \Delta u -  \frac{\partial u}{\partial t} =  -u \left( 2 \frac{\partial u}{\partial x} + \frac{\partial u}{\partial y} \right) +f(x,y), \\
u(x, y, t)=u(x+4, y, t),  \ u(x,-2,t)=-4, \ u(x,2,t)=2, \\
 u(x,y,0)=u_{init}(x,y,\mu ),\\
 x \in [-2, 2], \  y \in [-2, 2], \  t \in [0,1].
\end{cases}
\end{align}

We explicitly find the zero-order outer asymptotic functions:
\begin{align*}
&\varphi^{(-)} (x,y)= - \frac{2}{\sqrt{3 \pi }} \sqrt{
    \begin{aligned}
    &\sin \left(\frac{\pi  (x+ y)}{4} \right)-\sin \left(\frac{\pi  \left(x-2y-6\right)}{4} \right) \\
    &+3 \sin \left(\frac{\pi  ( x- y)}{4} \right)-3 \sin \left(\frac{\pi  \left(x-2y-2\right)}{4} \right)+12 \pi
    \end{aligned}
  },\\
&\varphi^{(+)} (x,y) = \frac{2}{\sqrt{3 \pi }} \sqrt{
    \begin{aligned}
    &\sin \left(\frac{\pi  ( x+ y)}{4} \right)-3 \sin \left(\frac{\pi  \left(x-2y+2\right)}{4} \right)\\
    &+3 \sin \left(\frac{\pi  ( x- y)}{4} \right)-\sin \left(\frac{ \pi  \left(x-2y+6\right)}{4}\right)+3 \pi
    \end{aligned}
  }.
\end{align*}

\begin{figure}[H]
\begin{center}
\includegraphics[width=0.4\linewidth,height=0.4\textwidth,keepaspectratio]{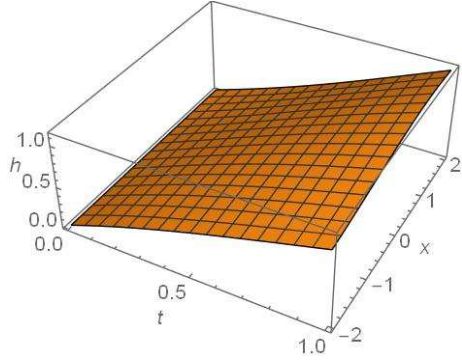}
\caption{Numerical solution to \eqref{h0mainequation} for $t \in [0,  1] $, $k=2$.}
\label{fig:x0example1}
\end{center}
\end{figure}

From the numerical solution to \eqref{h0mainequation}, we determine that the transition layer is located within the region $-2 \leq h_0 (x,t) \leq 2 $ for any $(x,t) \in [-2,2] \times [0,1]$ (see Fig. \ref{fig:x0example1}); thus, Assumption \ref{A3} is satisfied.

According to \eqref{InitialEx}, the initial function for problem \eqref{forwardexample1} which satisfies Assumption \ref{A4} takes the  form  $\displaystyle u_{init} (x,y,\mu )$ $\displaystyle =3\tanh\left(x+\frac{y}{0.08}\right)-1.$

\begin{figure}[H]
 \centering
\subfigure[]{
\includegraphics[width=0.4\linewidth]{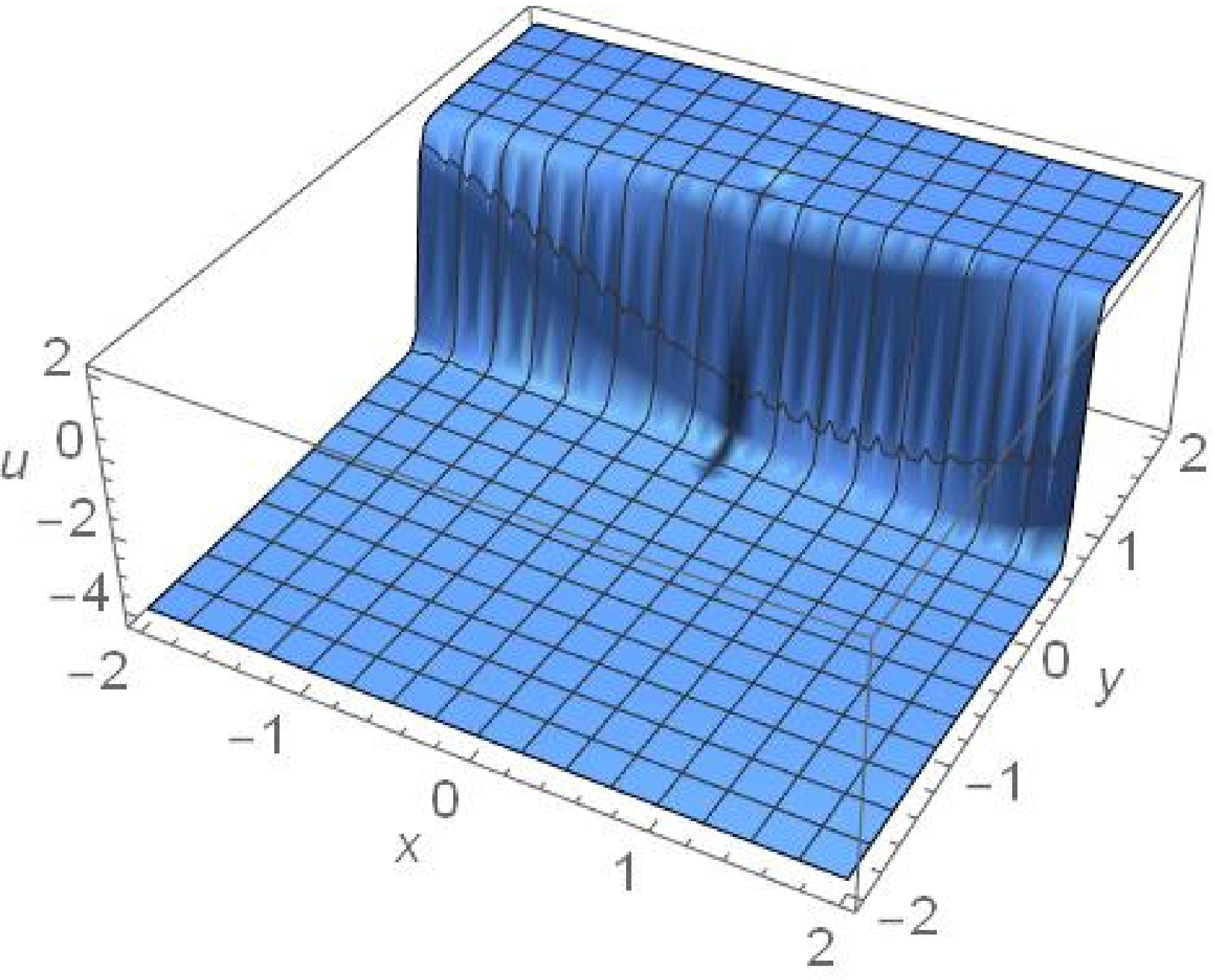} \label{fig:asymptoticsolexample1} }
\hspace{0ex}
\subfigure[]{
\includegraphics[width=0.4\linewidth]{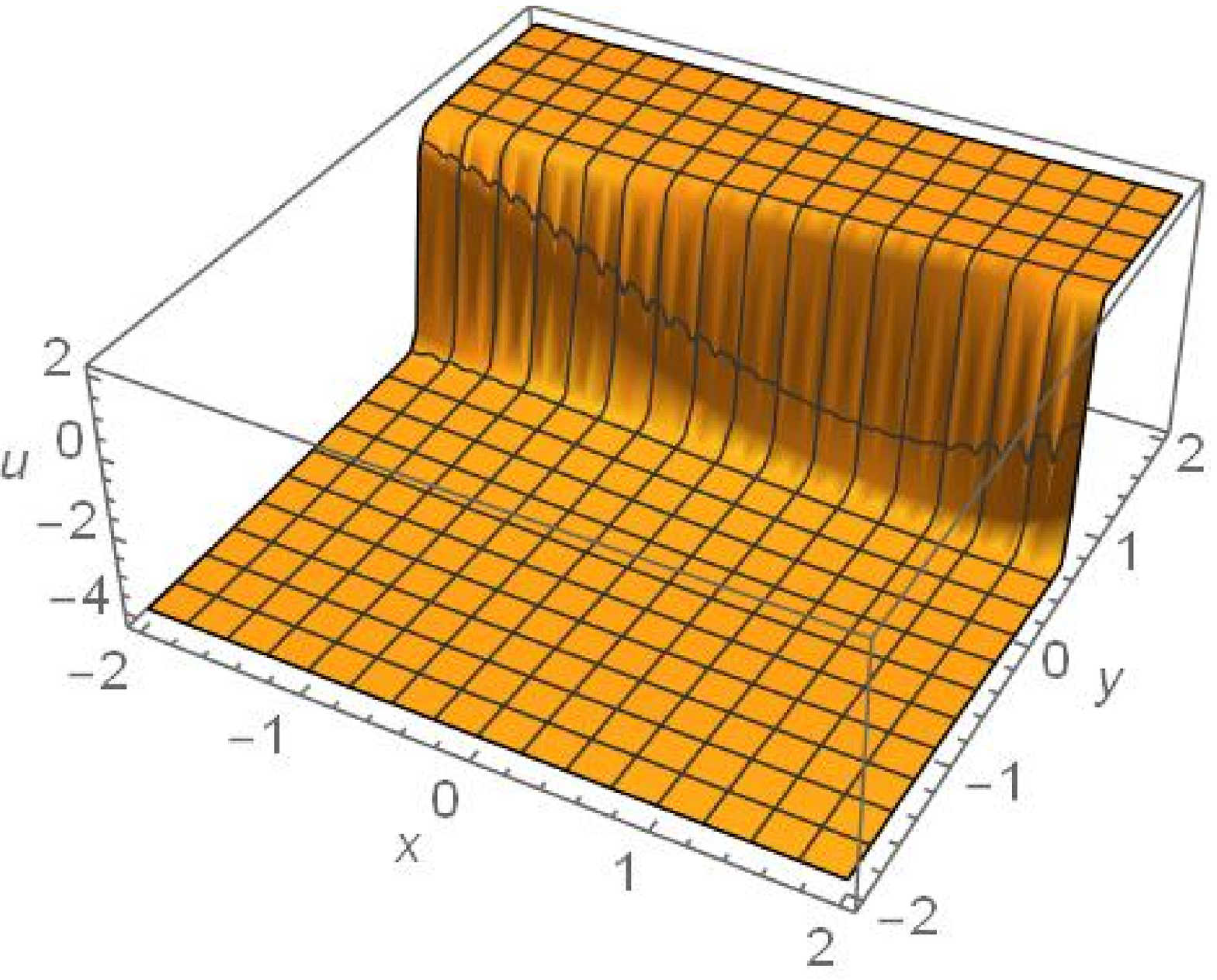} \label{fig:numericsolexample1} }
\caption{Asymptotic solution \subref{fig:asymptoticsolexample1} and numerical solution (using the finite-volume method) \subref{fig:numericsolexample1} of PDE \eqref{forwardexample1}  for $t_0=0.7$.} \label{solutionsexample1}
\end{figure}

Thus, all the assumptions formulated in this paper are satisfied, and the solution to considered equation \eqref{forwardexample1} is in the form of an autowave with a transitional moving layer localized near $h_0(x,t)$. The solution with the asymptotic expansion method for  $t \in [0,1]$ and $ x \in [-2,2]$ in the zero approximation has the form \eqref{asymptoticsolEXAMPLE1}, and is shown for $t_0=0.7$ in Fig. \ref{fig:asymptoticsolexample1}.

We also draw the numerical solution (using the finite-volume method) to problem \eqref{forwardexample1} at the moment $t_0=0.7$ on Fig. \ref{fig:numericsolexample1},  which we will use for the problem of identifying the source function.

The relative error of the asymptotic solution is $ \frac{\| U_0(x,y,t_0) - u(x,y,t_0) \|_{L^{2}([-2,2] \times [-2,2])}}{\|u(x,y,t_0) \|_{L^{2}([-2,2] \times [-2,2])}} = 0.0339. $

\subsubsection{Inverse problem}
We now consider  the  problem  of identifying the  source  function $f(x,y)$ in problem \eqref{forwardexample1}.

We use the following parameters: $t_0=0.7$,  $\delta = 1\%$, $n=50$, $m^{(-)} = 31$, $m^{(+)}= 39$, $m=50$.

\begin{figure}[H]
 \centering
\subfigure[]{
\includegraphics[width=0.4\linewidth]{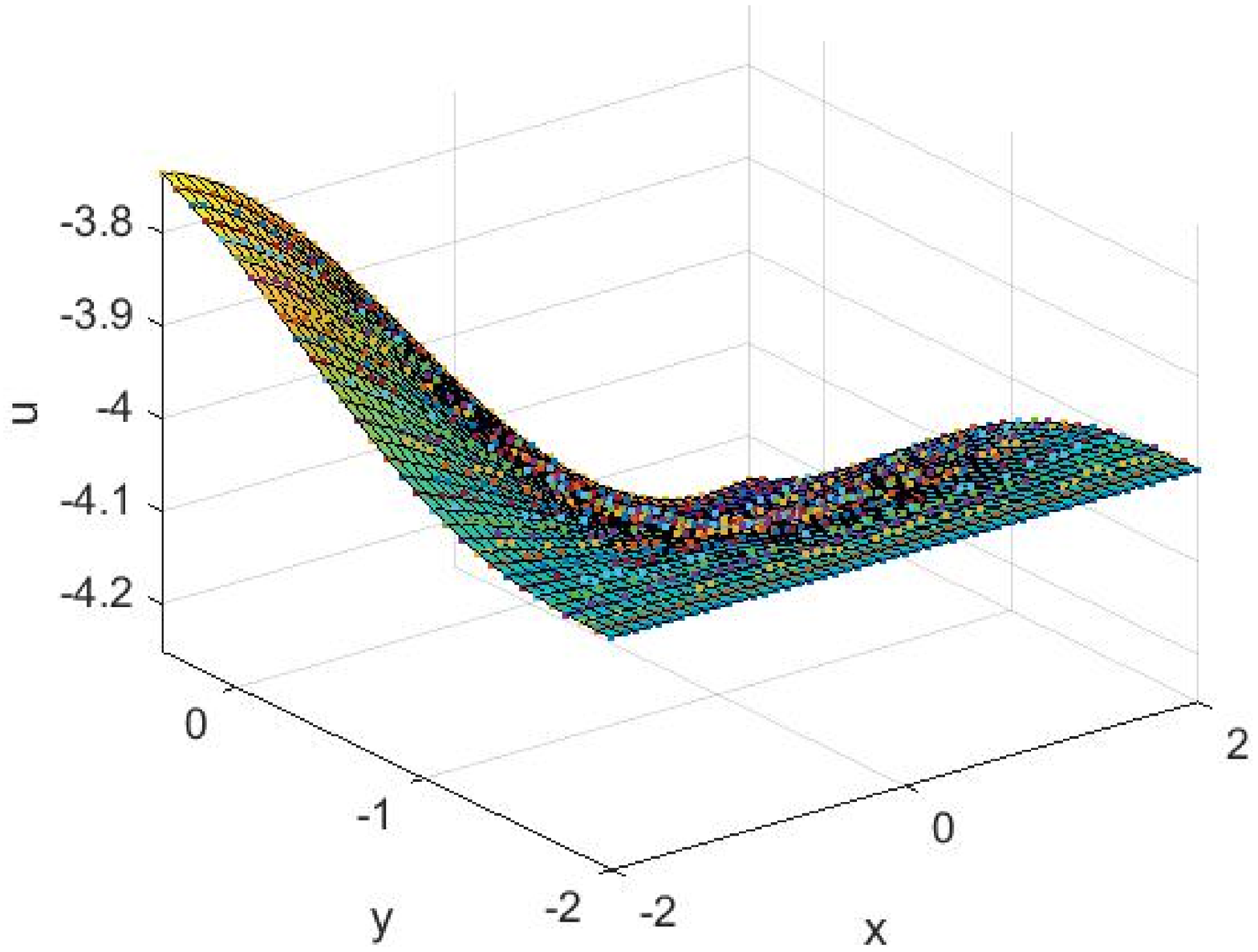} \label{fig:Lspline10-3} }
\hspace{0ex}
\subfigure[]{
\includegraphics[width=0.4\linewidth]{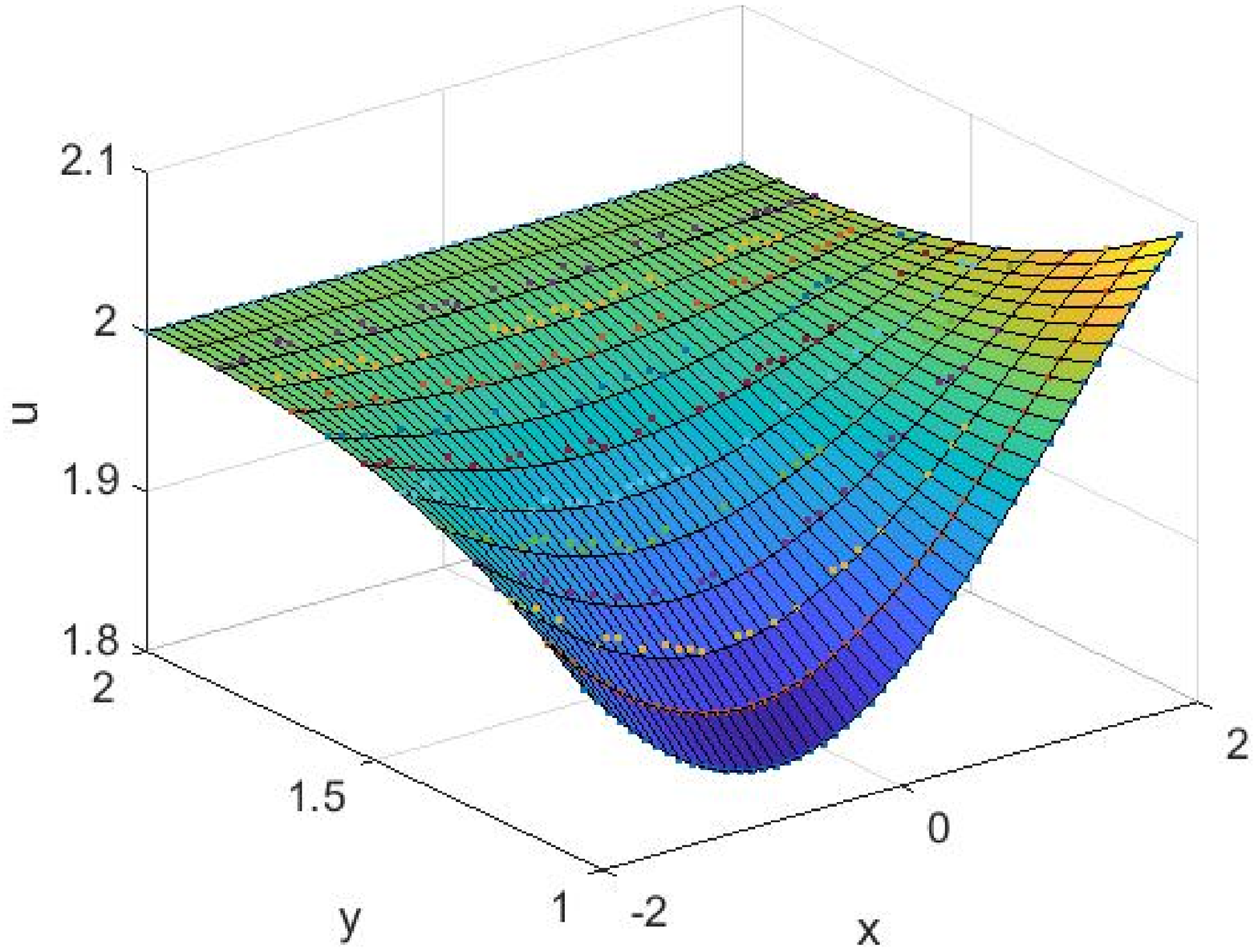} \label{fig:Rspline10-3} }
\caption{Left \subref{fig:Lspline10-3} and right \subref{fig:Rspline10-3} intervals, where surfaces represent smooth approximate data $u^\varepsilon(x,y,t_0)$, and points represent  noisy data   $ \{\{u^{\delta}_{i,j}(t_0)\}^{n}_{i=0} \}_{j=0}^{m^{(-)}}$, $\{ \{u^{\delta}_{i,j}(t_0)\}^{n}_{i=0} \}_{j=m^{(+)}}^{m}$ with $\delta=1 \%  $.} \label{fig:ualpha10-3Ex1}
\end{figure}

We take the data $u(x_i, y_i,t_0)$ obtained from the forward problem with the numerical method on the regions $[-2, h_0 (x,t_0)-\Delta h/2] $ and $[ h_0 (x,t_0)+\Delta h/2,2] $, and add random Gaussian noise with noise level $\delta$ to produce noisy data. Then, we smooth the obtained artificial noisy data on both intervals and obtain the smooth function $u^\varepsilon(x,y,t_0)$   (see Fig. \ref{fig:ualpha10-3Ex1}).

\begin{figure}[H]
\vspace{-1ex} \centering
\subfigure[]{
\includegraphics[width=0.4\linewidth]{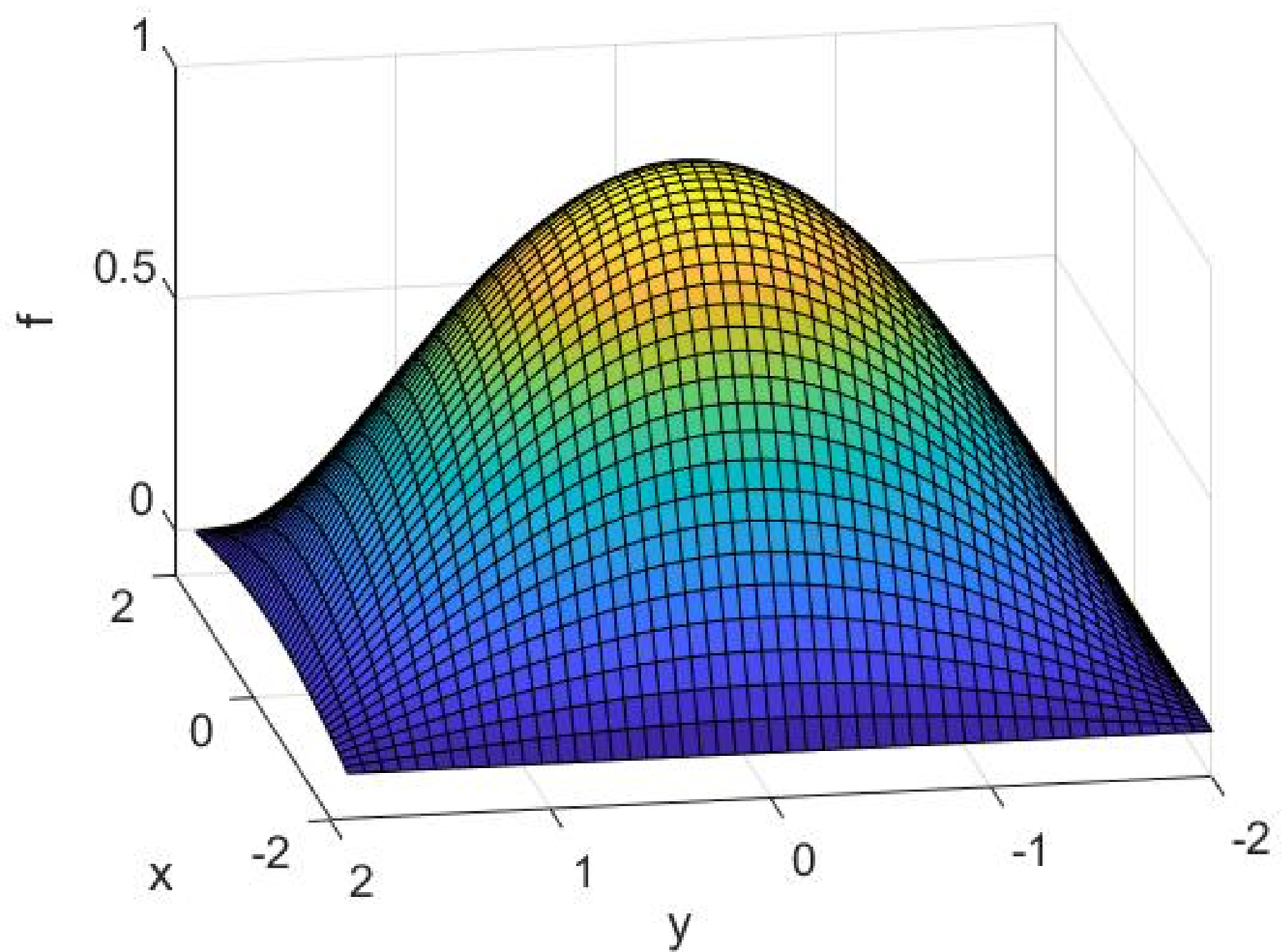} \label{fig:sourceorig} }
\hspace{0ex}
\subfigure[]{
\includegraphics[width=0.4\linewidth]{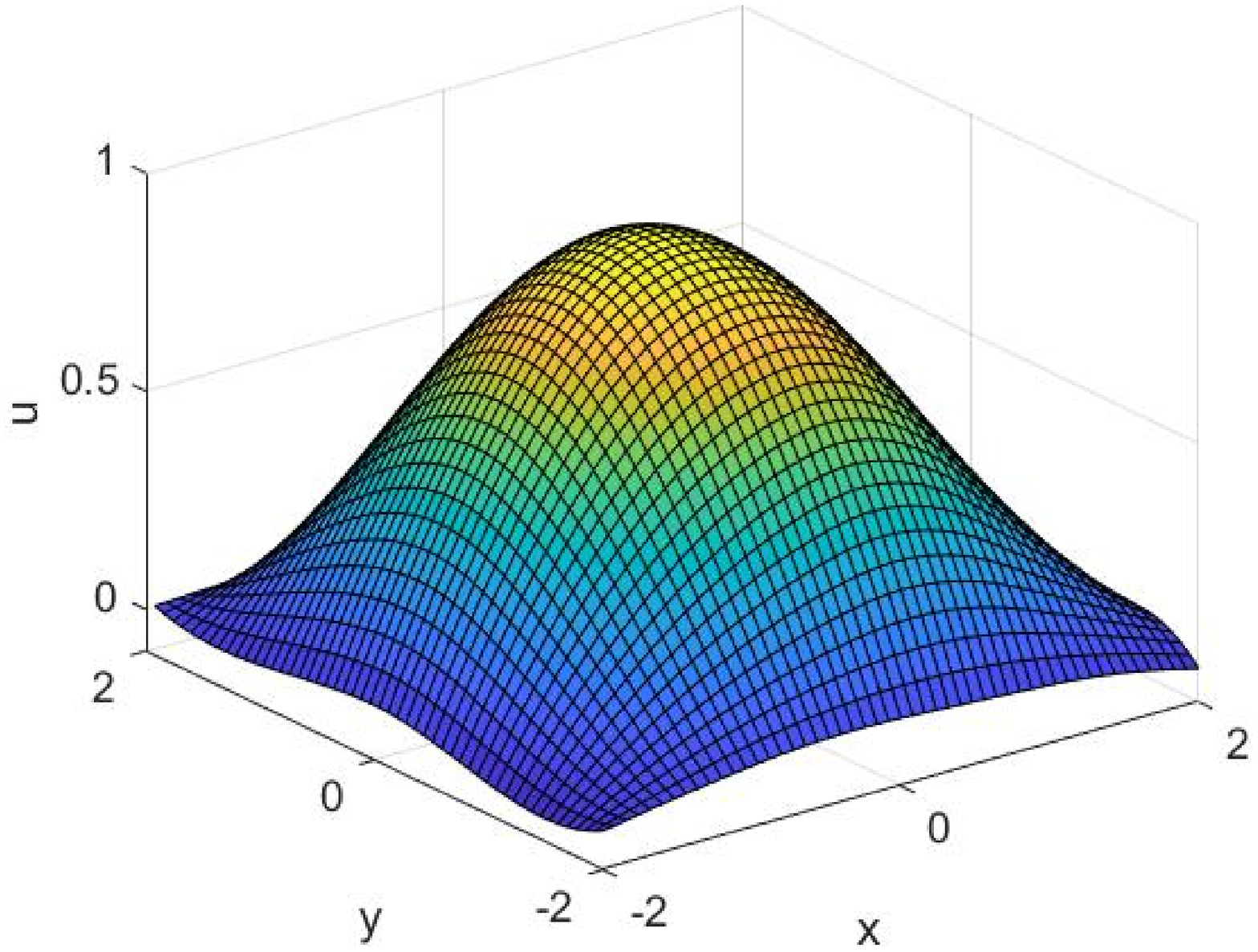} \label{fig:sourcereconstr} }
\caption{The result of reconstructing the source function $f^{\delta}(x,y)$ \subref{fig:sourcereconstr} for problem \eqref{forwardexample1}; this can be compared with the exact source function  $f^*=\cos{ ( \frac{\pi x}{4} )} \cos{( \frac{\pi y}{4} )}$ \subref{fig:sourceorig}.} \label{fig:sourcereconstruction}
\end{figure}

Following Algorithm \ref{alg:Framwork}, the regularized approximate source function $f^\delta(x,y)$ calculated using formula \eqref{fdelta} is drawn in Fig. \ref{fig:sourcereconstr}.

The relative error of the recovered source function is $ \frac{\| f^{\delta} - f^* \|_{L^2( [-2,2] \times [-2,2])}}{\|f^* \|_{L^2( [-2,2] \times [-2,2])}} = 0.0814$, from which we can conclude that our approach is stable and accurate.

\subsection{Example 2}
\subsubsection{Forward problem}

We consider equation \eqref{forwardexample1} with another source function $f(x,y)=y-2\cos{ (4 \pi x )} $ and a periodic boundary condition along the $ x $ axis:
\begin{align} \label{forwardexample2}
\begin{cases}
\displaystyle 0.08 \Delta u -  \frac{\partial u}{\partial t} =  -u \left( \frac{\partial u}{\partial x} + \frac{\partial u}{\partial y} \right) +f(x,y) ,\\
u(x, y, t)=u(x+2, y, t),  \  u(x,-1,t)=-8, \ u(x,1,t)=4,  \\
\displaystyle  u(x,y,0)=u_{init}(x,y,\mu ),\\
x \in [-1, 1], \  y \in [-1, 1], \  t \in [0,0.3].
\end{cases}
\end{align}

The zero-order asymptotic outer functions have the following form:
\begin{align*}
&\varphi^{(-)} (x,y)= -\frac{\sqrt{\sin (4 \pi  (x-y-1))-\sin (4 \pi  x)+\pi  y^2+63 \pi }}{\sqrt{\pi }},\\
&\varphi^{(+)} (x,y) = \frac{\sqrt{\sin (4 \pi  (x-y+1))-\sin (4 \pi  x)+\pi  y^2+15 \pi }}{\sqrt{\pi }}.
\end{align*}

Solving \eqref{h0mainequation}, we find the location of the transition layer for problem \eqref{forwardexample2}:
\begin{figure}[H]
\begin{center}
\includegraphics[width=0.4\linewidth,height=0.4\textwidth,keepaspectratio]{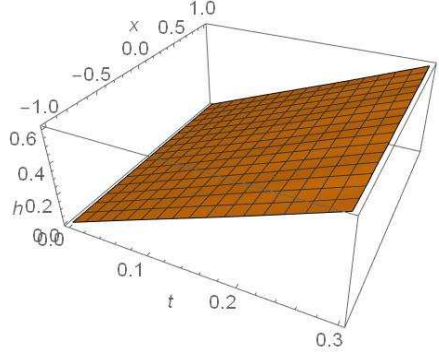}
\caption{Numerical solution of \eqref{h0mainequation} for $t \in [0,  0.3] $.}
\label{fig:x0example2}
\end{center}
\end{figure}
The transition layer is located within the region $-1 \leq h_0 (x,t) \leq 1 $ for any $(x,t) \in [-1,1] \times [0,1]$  (see Fig. \ref{fig:x0example2}), and thus Assumption \ref{A3} is satisfied.
The initial function for problem \eqref{forwardexample2}, which satisfies Assumption \ref{A4}, takes the form $\displaystyle u_{init} (x,y,\mu )=6\tanh\left(x+\frac{y}{0.08}\right)-2.$

\begin{figure}[H]
 \centering
\subfigure[]{
\includegraphics[width=0.4\linewidth]{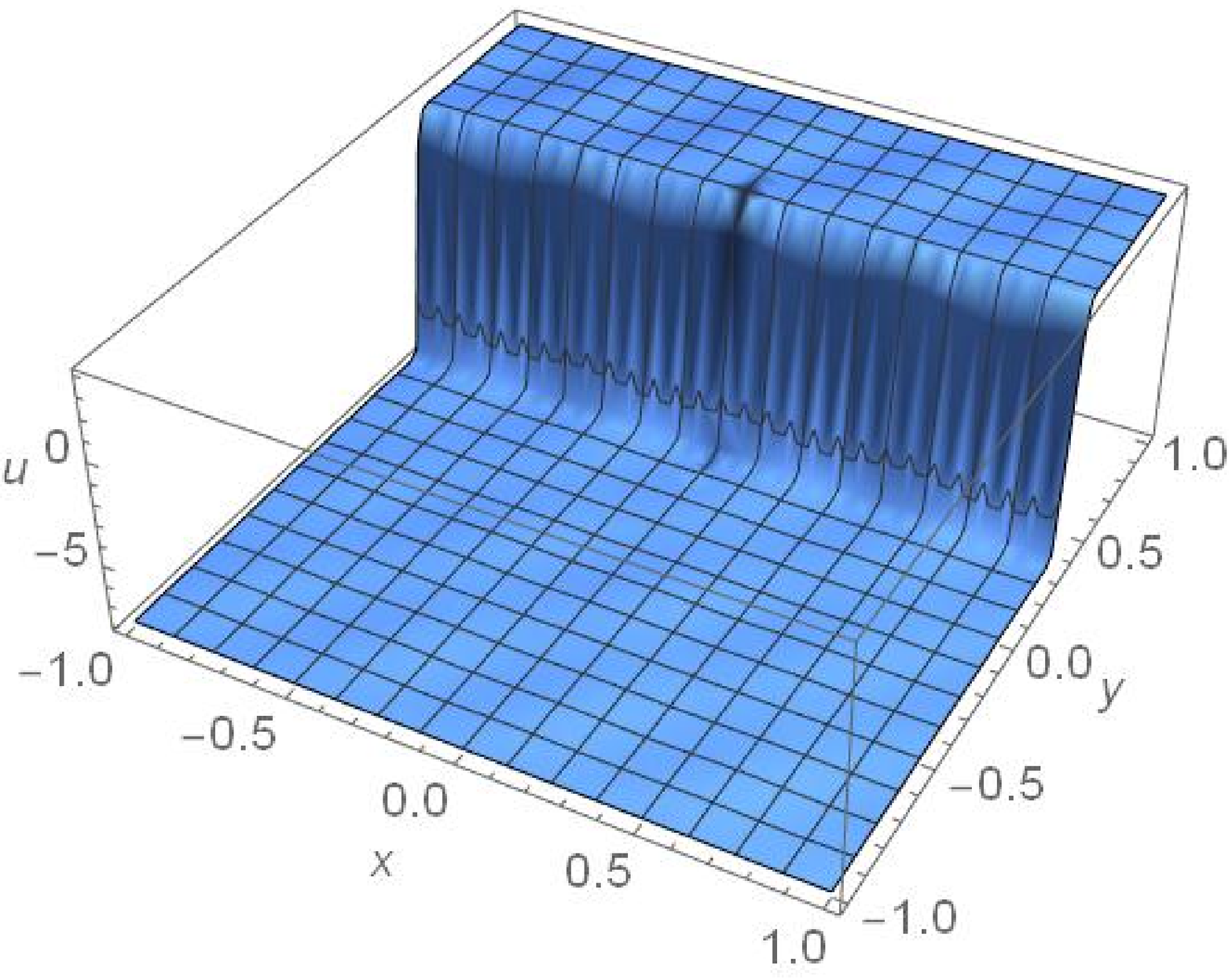} \label{fig:asymptoticsolexample2} }
\hspace{0ex}
\subfigure[]{
\includegraphics[width=0.4\linewidth]{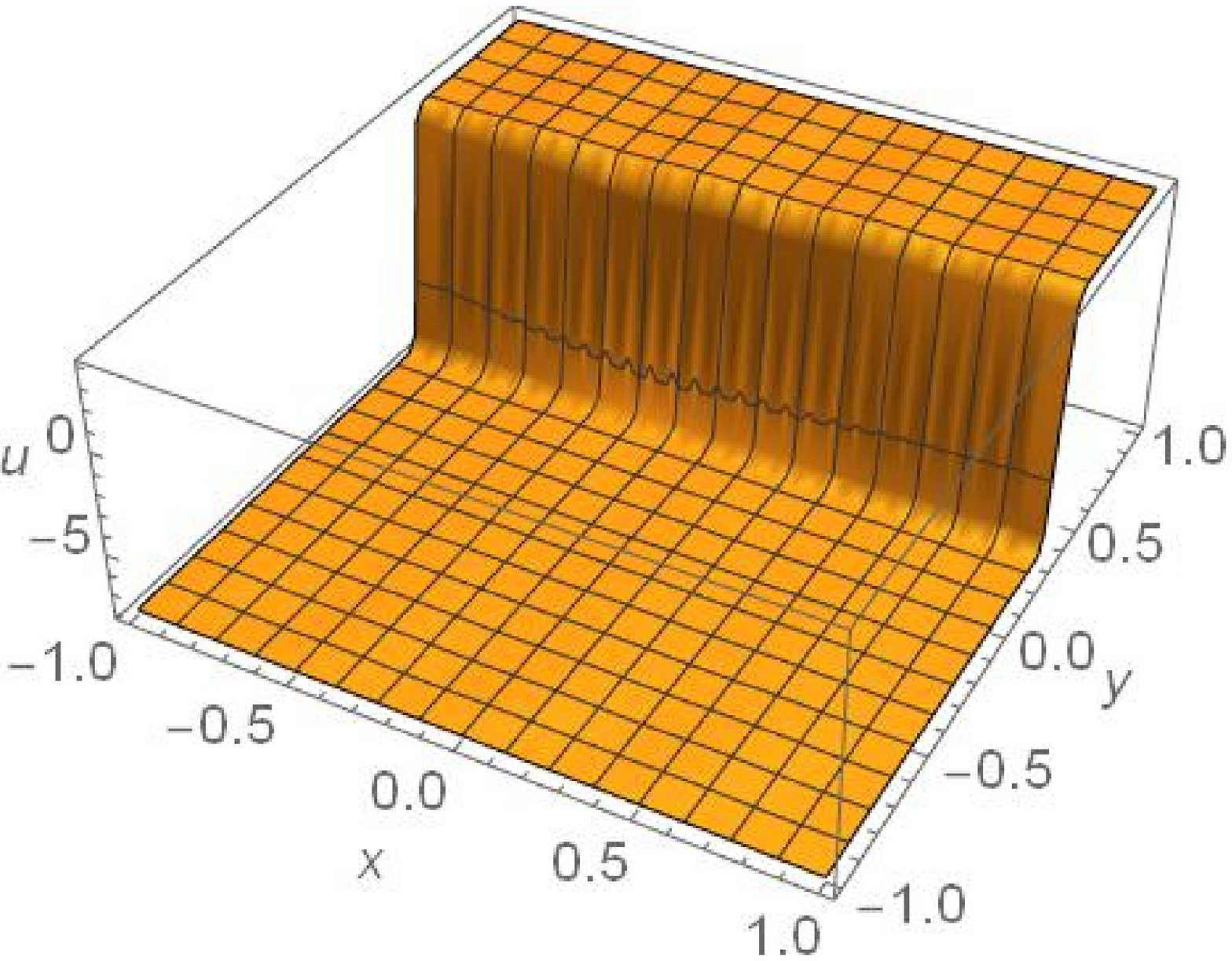} \label{fig:numericsolexample2} }
\caption{Asymptotic solution \subref{fig:asymptoticsolexample2} and numerical solution (using the finite-volume method) \subref{fig:numericsolexample2} of PDE \eqref{forwardexample2}  for $t_0=0.2$.}
\end{figure}

Thus, all the assumptions are satisfied and the considered equation, \eqref{forwardexample2}, has a solution in the form of an autowave; we represent the solution with the  asymptotic expansion method \eqref{asymptoticsolEXAMPLE1}   for $t_0=0.2$ in Fig. \ref{fig:asymptoticsolexample2} and the numerical solution in Fig. \ref{fig:numericsolexample2} (which we will use for the problem of identifying the source function).

The relative error of the asymptotic solution is  $\frac{\| U_0(x,y,t_0) - u(x,y,t_0) \|_{L^{2}([-1,1] \times [-1,1])}}{\|u(x,y,t_0) \|_{L^{2}([-1,1] \times [-1,1])}} = 0.0408 $.

\subsubsection{Inverse problem}
We now consider  the  problem  of  identifying the  source  function $f(x,y) $ of  problem \eqref{forwardexample2}. The following parameters are used in the simulation: $t_0=0.2$,  $\delta = 1\%$,  $n=50$, $m^{(-)} =28$, $m^{(+)}=34$, $m=50$.

\begin{figure}[H]
 \centering
\subfigure[]{
\includegraphics[width=0.4\linewidth]{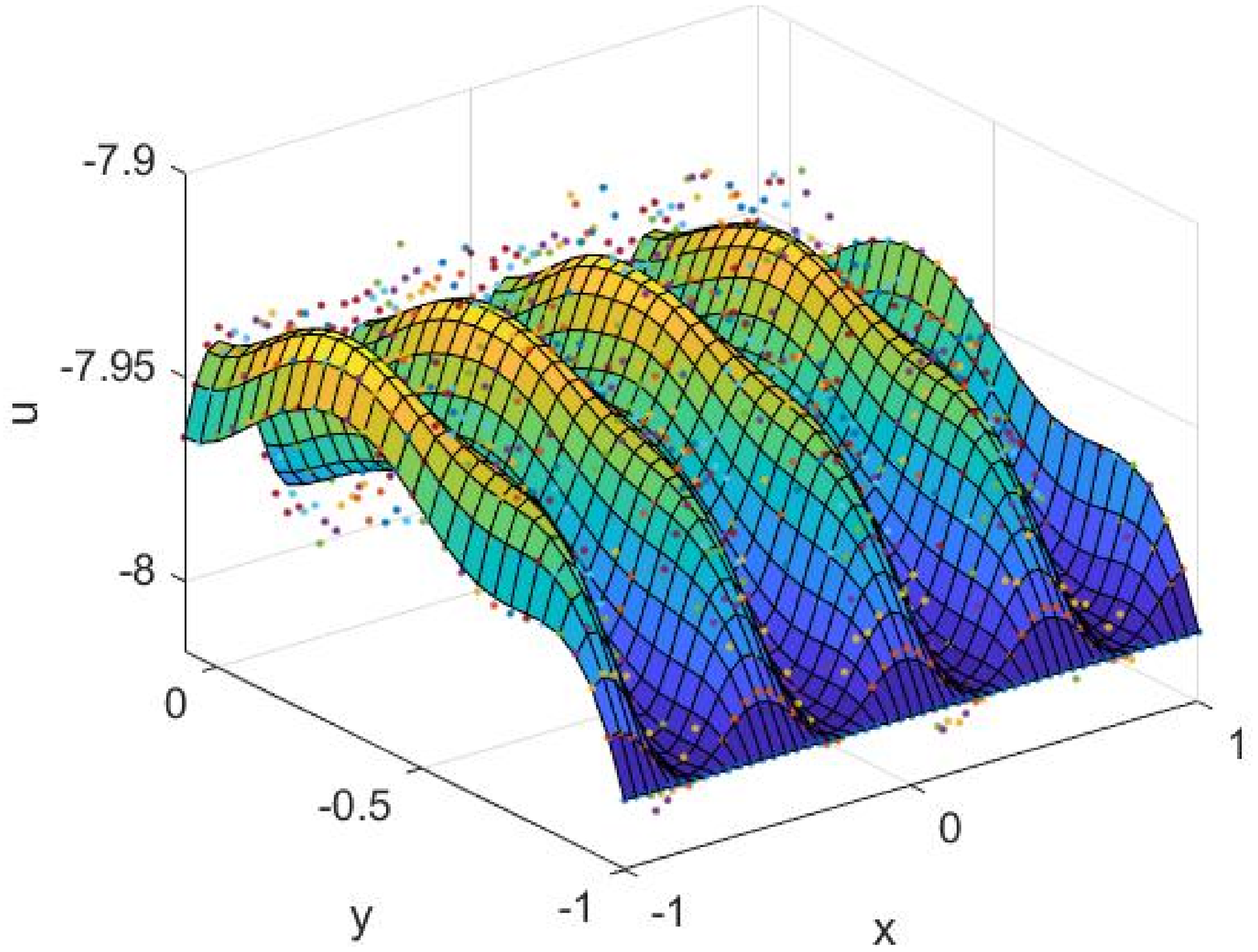} \label{fig:Lspline10-32} }
\hspace{0ex}
\subfigure[]{
\includegraphics[width=0.4\linewidth]{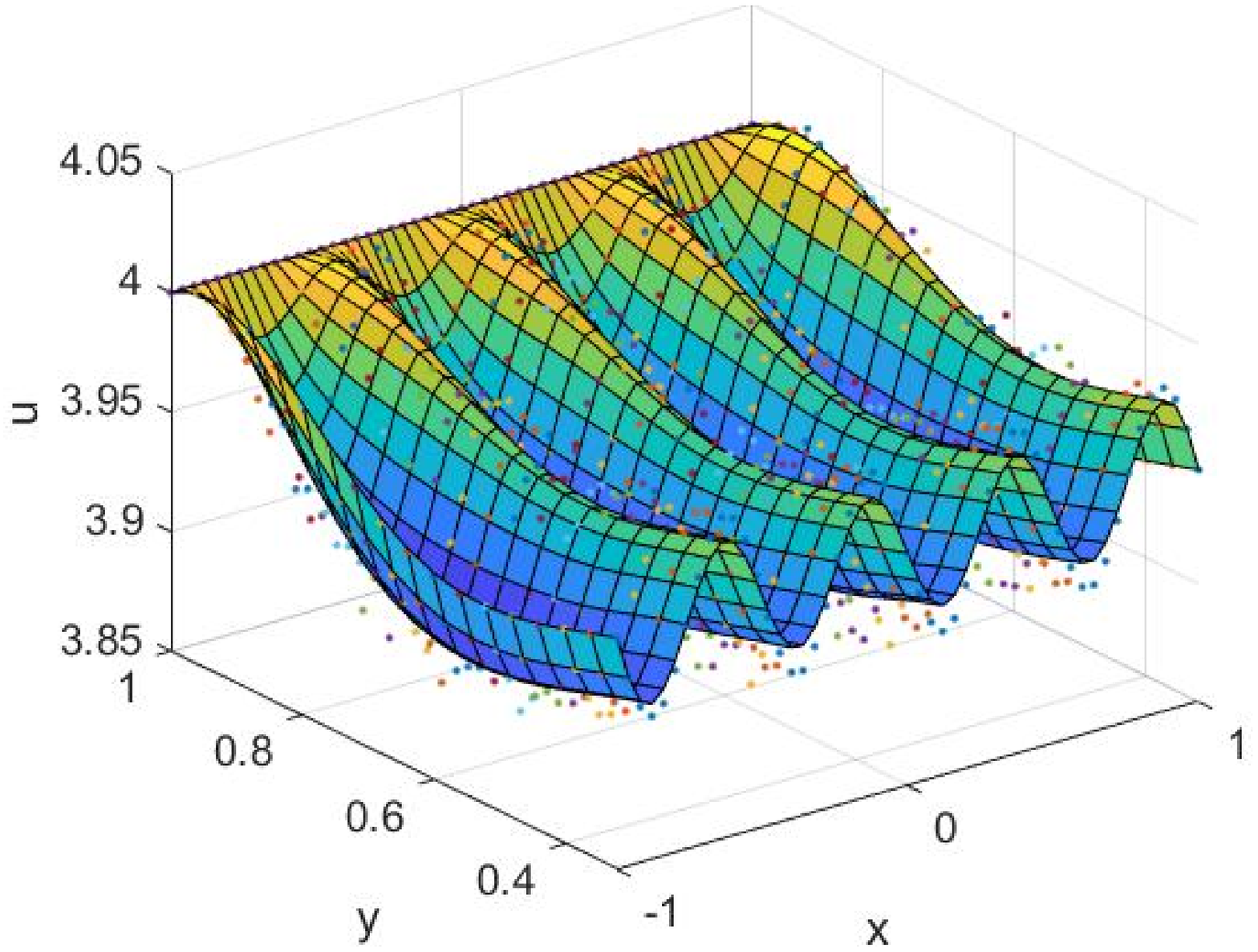} \label{fig:Rspline10-32} }
\caption{Left \subref{fig:Lspline10-32} and right \subref{fig:Rspline10-32} intervals where surfaces represent smooth approximate data $u^\varepsilon(x,y,t_0)$, and points represent noisy data   $ \{\{u^{\delta}_{i,j}(t_0)\}^{n}_{i=0} \}_{j=0}^{m^{(-)}}$, $\{ \{u^{\delta}_{i,j}(t_0)\}^{n}_{i=0} \}_{j=m^{(+)}}^{m}$, with $\delta=1 \%  $.} \label{fig:ualpha10-3}
\end{figure}

\begin{figure}[H]
\vspace{-1ex} \centering
\subfigure[]{
\includegraphics[width=0.4\linewidth]{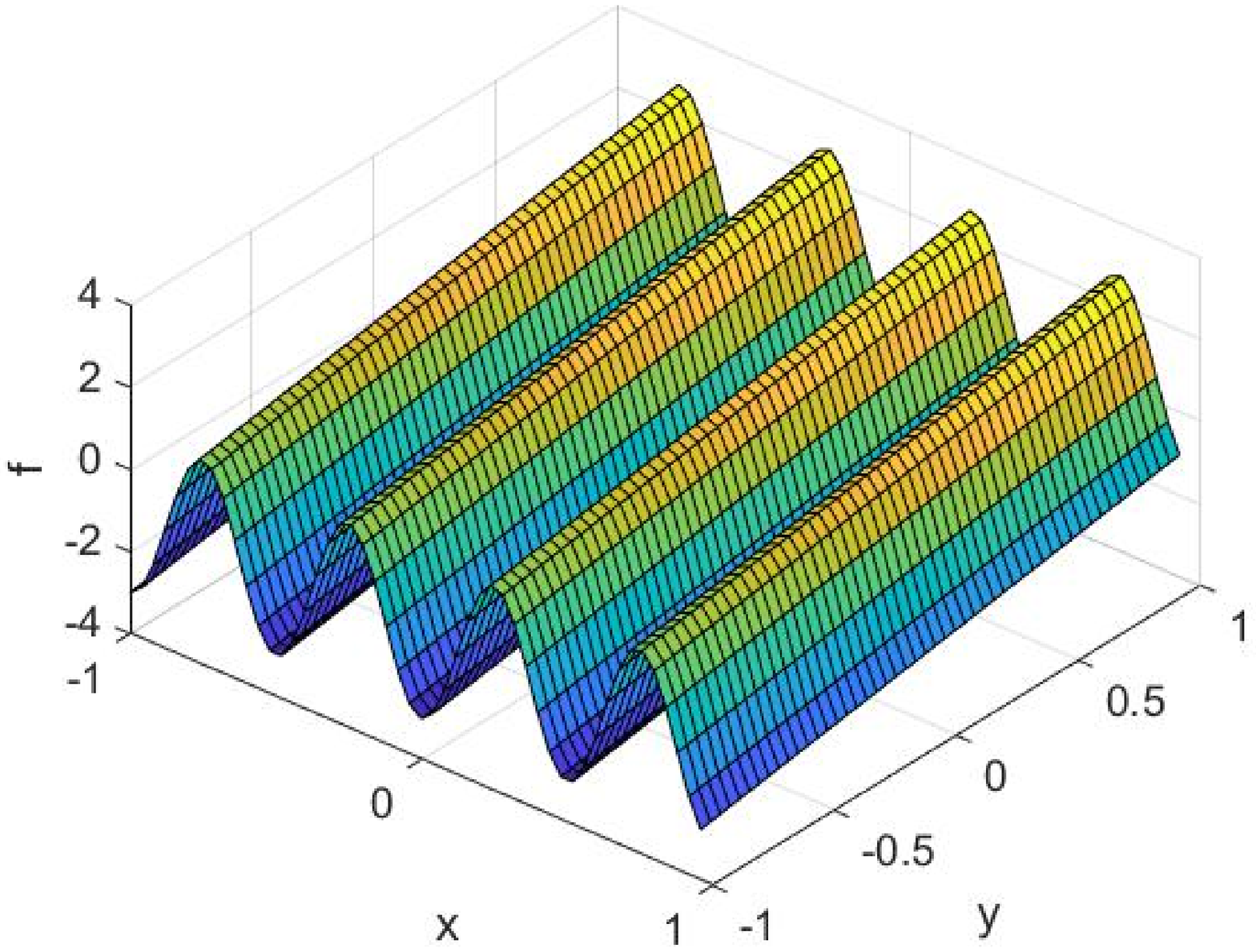} \label{fig:sourceorig2} }
\hspace{0ex}
\subfigure[]{
\includegraphics[width=0.4\linewidth]{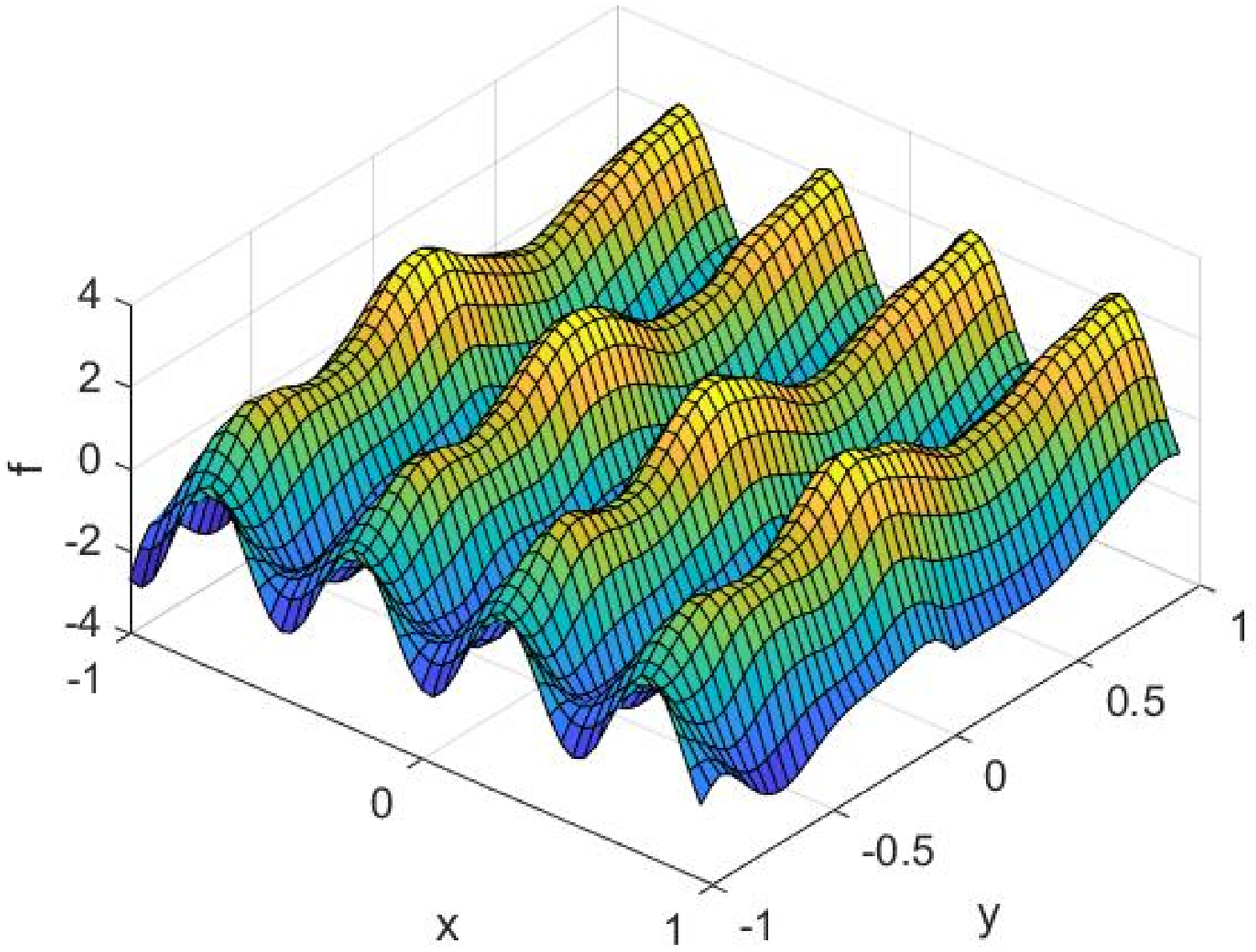} \label{fig:sourcereconstr2} }
\caption{The result of reconstructing the source function $f^{\delta}(x,y)$ \subref{fig:sourcereconstr2} for problem \eqref{forwardexample2}; this can be compared with exact source function $f^*=y-2\cos{ (4 \pi x )} $ \subref{fig:sourceorig2}.} \label{fig:sourcereconstruction2}
\end{figure}

Solving the inverse problem according to Algorithm \ref{alg:Framwork} for the artificial noisy data $u^{\delta}_{i,j}$ obtained from the forward problem \eqref{forwardexample2}, we obtain the smooth function $u^\varepsilon(x,y,t_0)$   (see Fig. \ref{fig:ualpha10-3}). According to Algorithm \ref{alg:Framwork}, the regularized approximate source function $f^\delta(x,y)$ is calculated using formula \eqref{fdelta}; it is drawn in Fig. \ref{fig:sourcereconstr2}.

The relative error of the recovered source function is  $\frac{\| f^{\delta} - f^* \|_{L^2([-1,1] \times [-1,1])}}{\|f^* \|_{L^2([-1,1] \times [-1,1])}} = 0.3768. $

\section{Conclusion} \label{Conclusion}

In this paper, we develop an asymptotic expansion-\hspace{0cm}regularization method for solving inverse source problems for two-dimensional time-dependent singularly perturbed PDEs. The essential idea is to use the asymptotic expansion methods to find a link equation between an unknown source function and the measurable quantities. By using the explored link equation, we avoid having to solve the PDE-constrained control problem. This simplification makes our inversion algorithm very efficient, and the methodology developed here can potentially be applied to a wide class of inverse problems in singularly perturbed PDEs.

%%%% Acknowledgments %%%%%%%%
\section*{Acknowledgments}
This work was funded by the National Natural Science Foundation of China (No. 12171036), Beijing Natural Science Foundation (Key project No. Z210001), Guangdong Fundamental and Applied Research Fund (No. 2019A1515110971) and Shenzhen National Science Foundation (No. 20200827173701001).

%%-----------------------------
%%      your bibliography
%%-----------------------------

\bibliographystyle{elsarticle-num}
\bibliography{references}

\end{document}